\documentclass[reqno]{amsproc}

\usepackage{amstext,amsmath,amssymb,amsfonts}
\usepackage[latin1]{inputenc}
\usepackage{epsfig}
\usepackage{hyperref}
\usepackage{color}

\usepackage{latexsym}

\newtheorem{lemma}{Lemma}
\newtheorem{corollary}{Corollary}

\newtheorem{definition}{Definition}
\newtheorem{theorem}{Theorem}

\newtheorem{proposition}{Proposition}

\newcommand{\vspi}{\vspace{0.4cm}}

\newcommand{\bea}{\begin{eqnarray}}
\newcommand{\eea}{\end{eqnarray}}
\newcommand{\beq}{\begin{equation}}
\newcommand{\eeq}{\end{equation}}
\newcommand{\enn}{\nonumber \end{equation}}

\newcommand{\tw}{ {\rm t}}

 \newcommand{\rk}{{\rm r}\,}

 \newcommand{\cG}{\mathcal{G}}
 \newcommand{\cV}{\mathcal{V}}
 \newcommand{\cE}{\mathcal{E}}
 \newcommand{\cF}{\mathcal{F}}

 \newcommand{\cR}{\mathcal{R}}
 \newcommand{\cP}{\mathcal{P}}

\newcommand{\bG}{\partial{\mathcal{G}}}
\newcommand{\bV}{{\mathcal{V}}_{\partial}}
\newcommand{\bE}{{\mathcal{E}}_{\partial}}
\newcommand{\bC}{{\mathcal{C}}_{\partial}}

 \newcommand{\mf}{\mathfrak{f}}
 \newcommand{\mh}{\mathfrak{h}}

 \newcommand{\sset}{\Subset}
 \newcommand{\eps}{\epsilon}

 \newcommand{\inter}{{\rm int}}
 \newcommand{\ext}{{\rm ext}}

\title[Recipe theorems for polynomial invariants on ribbon graphs]{Recipe theorems for polynomial invariants \\ on ribbon 
graphs with half-edges\footnote{ Preprint: ICMPA-MPA/2014/22}}

\author{Remi C. Avohou}
\address[R.C.A.]{
International Chair in Mathematical Physics and Applications,
ICMPA-UNESCO Chair, 072BP50, Cotonou, Rep. of Benin}
\email{remi.avohou@cipma.uac.bj}

\author{Joseph Ben Geloun}
\address[J.B.G.]{Max-Planck Institute for Gravitational Physics, 
Albert Einstein Institute, Am Muhlenberg 1, 14476, Potsdam,
Germany, and
International Chair in Mathematical Physics and Applications,
ICMPA-UNESCO Chair, 072BP50, Cotonou, Rep. of Benin}
\email{ jbengeloun@aei.mpg.de}

\author{Mahouton N. Hounkonnou}
\address[M.N.H.]{International Chair in Mathematical Physics and Applications,
ICMPA-UNESCO Chair, 072BP50, Cotonou, Rep. of Benin}
\email{norbert.hounkonnou@cipma.uac
.bj}

\begin{document}

\maketitle 

\begin{abstract}
We provide recipe theorems for the Bollob\`as and Riordan polynomial $\mathcal{R}$ defined on classes of ribbon graphs with half-edges introduced in arXiv:1310.3708[math.GT]. 
We also define a generalized transition polynomial $Q$ on this new category of ribbon graphs and  
establish a relationship between $Q$ and $\mathcal{R}$. 
\\

\end{abstract}

\tableofcontents

\section{Introduction}
The Bollob\`as-Riordan (BR) polynomial \cite{bollo} is a four-variable
polynomial  generalizing Tutte polynomial \cite{tutte} from simple graphs to ribbon graphs. Ribbon graphs are often called neighborhoods of graphs embedded into surfaces. Like Tutte polynomial, the BR polynomial satisfies a contraction/deletion recurrence relation 
and it is a universal invariant. The universality property of these invariants means that any invariant of graphs satisfying the same relations of contraction and deletion can be calculated from those. 
The universality can be also discussed in other contexts, for example, in statistical mechanics \cite{sok}
and even in quantum field theory \cite{Duchamp:2013pha,Duchamp:2013joa,Tanasa:2012pm}. 

In  \cite{joan}, the authors provide a ``recipe theorem'' for the BR polynomial intimately close to its universality property \cite{Bry,joan1,joan2}. This theorem is a very useful tool in order to evaluate any function $F$ on these ribbon graphs in terms of the BR polynomial itself, once one imposes that $F$ satisfies the same contraction/deletion recurrence relation and few more properties of the BR polynomial. Both proofs of the universality property and the recipe theorem of the BR polynomial are  based on the property of the contraction/deletion relation and the understanding of other ingredients such as chord diagrams associated with one-vertex ribbon graphs.

In a subsequent work \cite{rca}, another BR polynomial $\cR$ is introduced on a new class of ribbon graphs called half-edged ribbon graphs (HERGs). A half-edge or half-ribbon (HR) is  a ribbon incident to a unique vertex without forming a loop. The presence of HRs in a ribbon graph have 
notable combinatorial properties which makes the polynomial 
found in \cite{rca} a nontrivial one. Furthermore, the presence HRs also allows one to define the cut operation of a ribbon edge which differs from the usual edge deletion. The new polynomial $\cR$ found on HERGs then satisfies a contraction/cut recurrence relation. There exists another interesting operation on ribbon graphs which consists in moving the 
HRs on the boundary of these graphs \cite{rca2}. It has been
proved that one can quotient the action of these HR
moves  and get the so-called HR-equivalent 
classes of ribbon graphs with half-ribbons. There exists
a natural extension of $\cR$ on these equivalence
classes. The main result in \cite{rca2} is the proof of 
the universal property of the polynomial $\cR$ on HERGs
and of its extension to HR-equivalent classes
of ribbon graphs. This statement
 relies on the understanding and generalization of the tools 
necessary to  the proof of the universality of the original
BR polynomial. 
 
 The polynomial $\cR$ is universal on HERGs or on HR-equivalent 
ribbon graphs (in the following discussion, we will call these HR-classes
or simply classes). HR-classes are in a sense more fundamental.
 Indeed, as far as one is concerned with the evaluation of $\cR$ on a HERG $\cG$, we have $\cR(\cG)= \cR([\cG])$ where
$[\cG]$ is the HR-class of $\cG$. 
Hence, we might not need to encode all the information 
about the positions of the HRs on a given ribbon graph
before evaluating its invariant $\cR$. 
We can address the next question related to the existence
of a universal property of $\cR$: ``Can one provide a recipe
theorem for $\cR$ on these HR-classes of ribbon graphs?''. Answering this question is one of the main purposes of this paper.

In the present work which should be considered as a companion paper of \cite{rca2}, we will provide recipe theorems (Theorems \ref{theo:recipe} and \ref{theo:recipe2}) for computing any function $F$ on 
HR-classes satisfying the contraction/cut rule from the knowledge of $\cR$. 
Our proof is then  different from the one introduced
in \cite{joan}. Indeed, the authors of this contribution
based their proof on four items among which (item 2 therein) the factorization property of the BR polynomial
when evaluated on one-point-joint ribbon graphs. To be clearer, 
consider $G_1$ and $G_2$ two distinct ribbon graphs (without any half-edge consideration) $R(G_1 \cdot G_2) = R(G_1)R(G_2)$, where $R$ is the BR polynomial in the original sense of \cite{bollo} and $G_1 \cdot G_2$ is the one-point-joint operation of $G_1$ and $G_2$. With this property,
the way to evaluate $R$ on particular chord diagrams becomes simple. 
It turns out that 
the polynomial $\cR$ defined on half-edged ribbon graphs does not satisfy the same property, namely,
$\cR(G_1 \cdot G_2)\neq \cR(G_1)\cR(G_2)$, for
 $G_1$ and $G_2$ ribbon graphs with half-edges. This makes the polynomial $\cR$ a radically different invariant. 
We have identified new sets of conditions (replacing in particular
the failing item 2) under which a new genuine recipe can be provided.  In fact, one of the recipes found here can be considered as truly fundamental in the following
sense: our method yielding Theorem \ref{theo:recipe2} 
can be slightly adjusted for the polynomial $R$ and we will obtain the recipe as determined in \cite{joan}.  Finally, as a second 
main result of this work, 
we introduce a generalized version of the
transition polynomial \cite{joan,joan1} for HERGs and establish another
main result  (Theorem \ref{Theo:transpoly}) which provides 
a relationship between the transition polynomial found
and $\cR$.  The extension of this invariant
to HR-classes is immediate.

This paper is organized as follows. In section \ref{sect:hrequiv}, we recall some results on the BR polynomial for HERGs and for HR-classes and its universality property obtained in \cite{rca2}. 
The reader must be aware of the basics of ribbon graphs 
as found in \cite{bollo} or \cite{joan1} (but for notations
closer to the present paper, we refer to \cite{rca2}). 
In section \ref{sect:Recipe}, we give the first main result which is the proof of recipe theorems for this polynomial. We finally define, in section \ref{sect:trans}, the generalized transition polynomial $Q$ on HERGs via the medial graph construction (associated with HERGs) and find a relationship between $Q$ and the BR polynomial $\cR$.

\section{Polynomial invariants}  
\label{sect:hrequiv}

This section recalls the main results of \cite{rca2}: 

- the definition of ribbon graphs with half-edges
(studied originally in \cite{krf}) then precise the HR-equivalent classes that we will consider in this work; 

- the definition of a polynomial invariant 
on HR-classes of ribbon graphs and its universality property.

\begin{figure}[h]
 \centering
     \begin{minipage}[t]{.8\textwidth}
      \centering
\includegraphics[angle=0, width=3cm, height=1.5cm]{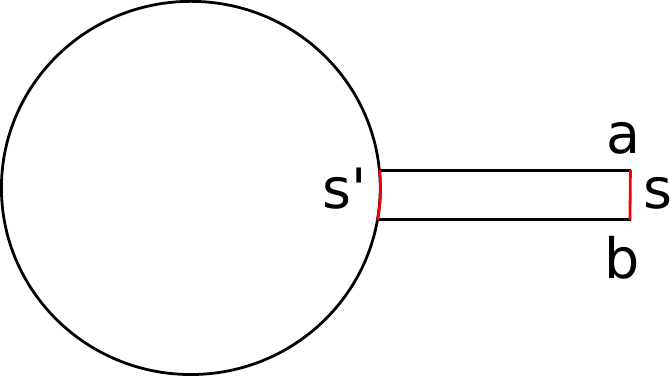}
\caption{ {\small A HR with two end segments (in red): $s'$
touching the vertex and $s$ external; the ends $a$ and $b$ of $s$
are the external points. }}
\label{fig:flag}
\end{minipage}
\end{figure}

\begin{definition}[Half-ribbon edges \cite{rca}]
\label{ribflag}
A half-ribbon edge (or simply half-ribbon, denoted henceforth HR) is a ribbon incident to a unique vertex by a unique segment and without forming loops. A  HR has two segments one touching a vertex 
and another free or external segment. The end-points of any free 
segment are called external points of the HR (see Figure \ref{fig:flag}).
\end{definition}

\begin{figure}[h]
\centering
\begin{minipage}[t]{.8\textwidth}
\centering
\includegraphics[angle=0, width=7cm, height=1cm]{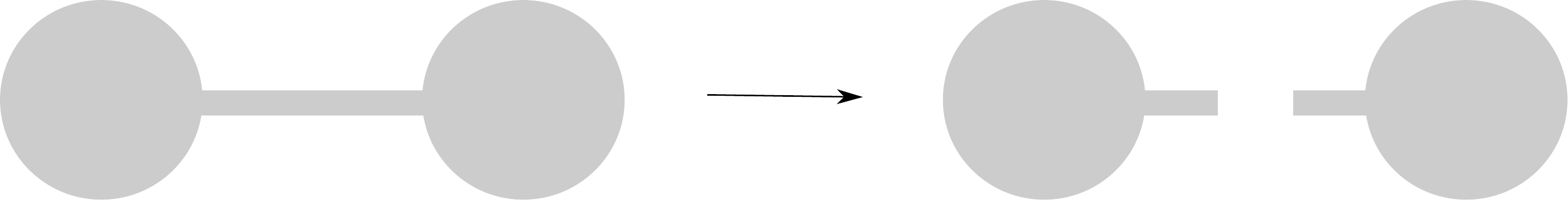}
\caption{ {\small Cutting a ribbon edge. }}
\label{fig:cut}
\end{minipage}
\end{figure}

\begin{definition}[Cut of a ribbon edge \cite{krf}]
Let $\cG$ be a ribbon graph and $e$ be an edge in $\cG$.
The cut graph $\cG \vee e$ is the graph obtained by 
removing $e$ and let two HRs attached at the end vertices
of $e$. If $e$ is a self-loop, the two HRs are on the same 
vertex. (See an illustration in Figure \ref{fig:cut}.)
\end{definition}

In a ribbon graph, we can identify three different kinds of edges: the bridge, the loop and the regular edge. Consider an edge $e$ of a ribbon graph $\cG$, $e$ is called a bridge in $\cG$ if its removal disconnects a component of $\cG$. If the two ends of $e$ are incident to the same vertex $v$ of $\cG$, $e$ is called a loop in $\cG$. A loop $e$ is trivial if there is no cycle in $\cG$ which can be contracted to form a loop $f$ interlaced with $e$. The edge $e$ is a regular edge of $\cG$ if it is neither a bridge nor a loop.

\begin{figure}[h]
\centering
\begin{minipage}[t]{.8\textwidth}
\centering
\includegraphics[angle=0, width=7cm, height=1.4cm]{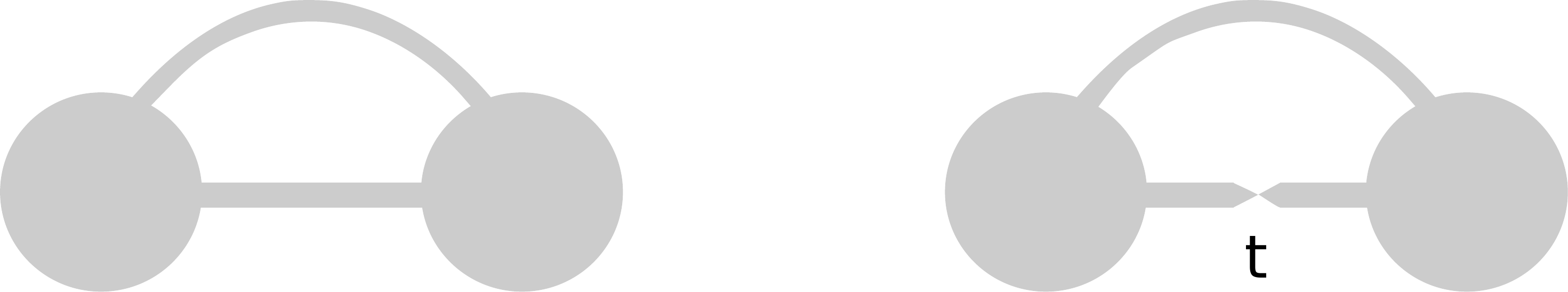}
\caption{ {\small Untwisted (left) and twisted (right) edge notations.
}}
\label{fig:edges}
\end{minipage}
\end{figure}

There are twisted and untwisted  ribbon edges (see illustration in Figure \ref{fig:edges}). A loop $e$ attached to a vertex $v$ of a ribbon graph $\cG$ is twisted if $v\cup e$ forms a M\"obius band as opposed to an annulus (an untwisted loop).

The notion of a half-edged ribbon graph (HERG) may be
now introduced.

\begin{figure}[h]
 \centering
     \begin{minipage}[t]{.8\textwidth}
      \centering
\includegraphics[angle=0, width=6.5cm, height=2cm]{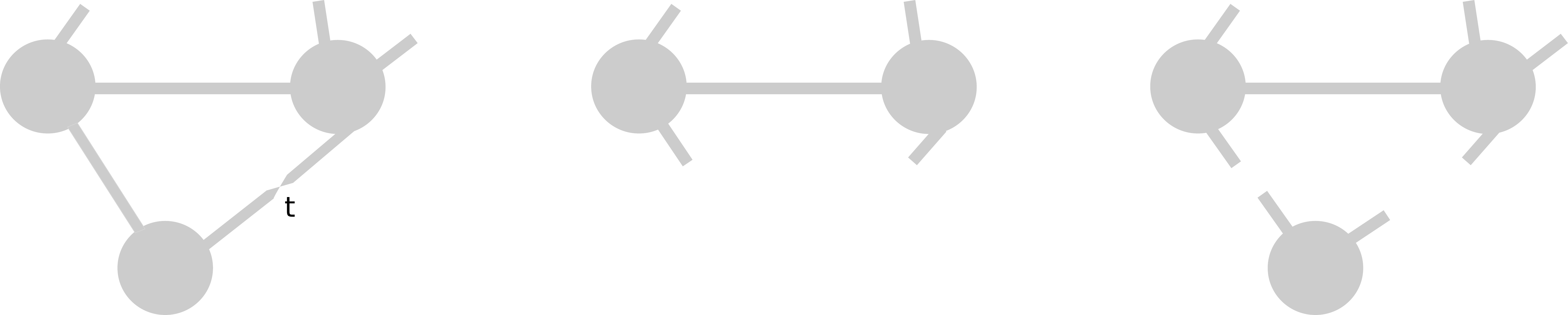}
\vspace{0.3cm}
\caption{ {\small A ribbon graph with HRs $\cG$, 
a c-subgraph $A$ and a spanning c-subgraph $\tilde {A}$. }}
\label{fig:ribgraph}
\end{minipage}
\put(-235,-10){$\cG$}
\put(-175,-10){$A$}
\put(-97,-10){$\tilde {A}$}
\end{figure}

\begin{definition}[Ribbon graph with HRs \cite{rca}]

$\bullet$ A ribbon graph $\cG$ with HRs  is a ribbon graph $\cG(\cV,\cE)$ with a set $\mf$ of HRs defined by the disjoint union 
of  $\mf^1$ the set of HRs obtained only from the cut of all edges of $\cG$ and a set $\mf^0$  of additional HRs together with a relation which 
associates  with each additional HR a unique vertex. 
We denote a ribbon graph with  set $\mf^0$ of additional HRs  as
$\cG(\cV,\cE, \mf^0)$. (See Figure \ref{fig:ribgraph}.)

$\bullet$ A c-subgraph $A$ of $\cG(\cV,\cE, \mf^0)$ is defined as a ribbon graph with HRs $A(\cV_A,\cE_A,\mf^0_A)$ the vertex 
set of which is a subset of $\cV$, the edge set of which  is a subset of $\cE$ together with their end vertices. 
Call $\cE_A'$ the set of edges incident to the vertices of $A$ and not contained in $\cE_A$.
The HR set of $A$  contains a subset of $\mf^0$ plus additional HRs attached to the vertices 
of $A$ obtained by cutting all  edges in $\cE_A'$. In symbols, $\cE_A \subseteq \cE$ and
$\cV_A \subseteq \cV$, $\mf^0_A = \mf^{0;0}_A \cup \mf^{0;1}_A(\cE_A)$
with $\mf^{0;0}_A \subseteq \mf^0$ and $\mf^{0;1}_A (\cE_A)\subseteq \mf^1$, where $\mf^{0;1}_A (\cE_A)$
 is the set of HRs obtained by cutting all edges in $\cE_A'$ 
and incident to vertices of $A$. 
We write $A \subseteq \cG$. 
 (See a c-subgraph $A$ illustrated in Figure \ref{fig:ribgraph}.)

$\bullet$ A spanning c-subgraph $A$ of $\cG(\cV,\cE,\mf^0)$ is defined as   a c-subgraph $A(\cV_A,\cE_A,\mf^0_A)$ of $\cG$ 
with all vertices  and all additional HRs of $\cG$. Hence $\cE_A\subseteq \cE$ and $\cV_A = \cV$, $\mf^0_A = \mf^{0} \cup \mf^{0;1}_A(\cE_A)$. 
 We write $A \sset \cG$. 
(See $\tilde A$ in Figure \ref{fig:ribgraph}.)
\end{definition} 

 A cutting-spanning subgraph can be obtained as follows: 
cut a subset of edges of a given graph. Then consider the spanning subgraph formed by the resulting graph. The set of HRs of this subgraph is the disjoint union of the set of HRs of the initial graph ($\mf^{0}$) plus an additional set induced by the cut of the edges.

Cutting an edge of a HERG brings some modifications on the boundary faces of this graph. We  obtain new boundary faces (following the contour of the HRs) which are  different from the ones which follow 
only the boundary of ribbon edges.

\begin{definition}[Closed and open faces \cite{Gurau:2009tz}]
\label{pinch}
Consider $\cG(\cV,\cE,\mf^0)$ a ribbon graph with HRs.

$\bullet$ A closed or internal face is a boundary face component of a ribbon graph (regarded as a geometric ribbon) 
which never passes through any free segment of additional HRs.  The set of closed faces is denoted $\cF_{\inter}$. 

$\bullet$ An open  or external face is a boundary face component 
leaving an external point of some HR rejoining another external point.  
 The set of open faces is denoted  $\cF_{\ext}$. 

$\bullet$ The two boundary lines of a ribbon edge or a HR
are called strands. Each strand belongs either to a closed or to open face. 

$\bullet$ The set of faces $\cF$ of a graph is defined by 
$\cF_{\inter} \cup \cF_{\ext}$.

$\bullet$ A graph is said to be open if $\cF_{\ext}\neq \emptyset$
i.e. $\mf^0\neq \emptyset$. It is closed otherwise. 

\end{definition}

\begin{figure}[h]
 \centering
     \begin{minipage}[t]{.8\textwidth}
      \centering
\includegraphics[angle=0, width=4cm, height=2cm]{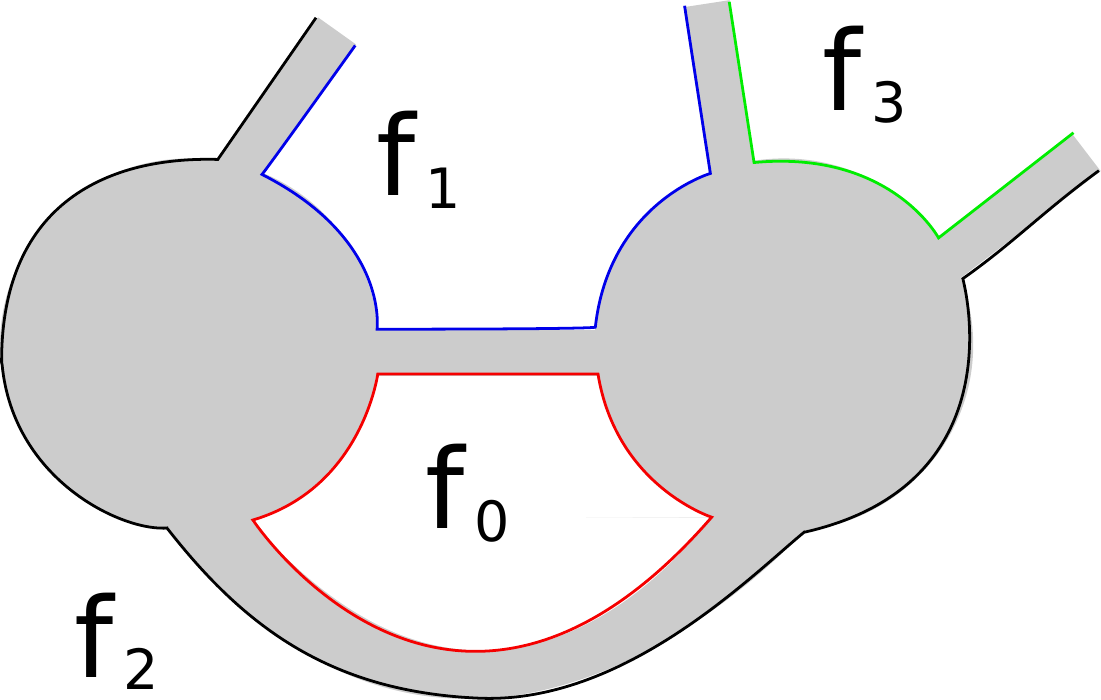}
\vspace{0.3cm}
\caption{ {\small A ribbon graph with set 
of internal faces 
$\cF_{\inter}=\{f_0\}$, and set of external 
faces $\cF_{\ext}=\{f_1,f_2,f_3\}$. }}
\label{fig:faces}
\end{minipage}
\end{figure}

Open and closed faces are illustrated in Figure \ref{fig:faces}.

\begin{definition}[Boundary graph \cite{Gurau:2009tz}]
\label{bnd}
$\bullet$ The boundary $\bG$ of a ribbon graph $\cG(\cV,\cE,\mf^0)$ 
is a simple graph $\bG(\bV,\bE)$
such that $\bV$ is one-to-one with $\mf^0$
and $\bE$ is one-to-one with $\cF_{\ext}$. 

$\bullet$ The boundary graph of a closed graph is empty. 
\end{definition}

\begin{figure}[h]
 \centering
     \begin{minipage}[t]{.8\textwidth}
      \centering
\includegraphics[angle=0, width=2cm, height=1cm]{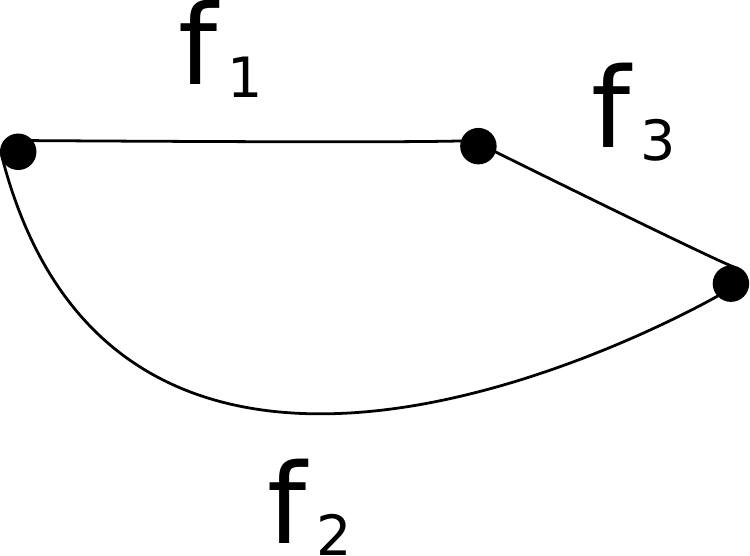}
\vspace{0.3cm}
\caption{ {\small The boundary graph associated with 
the ribbon graph in Figure \ref{fig:faces}. }}
\label{fig:boundary}
\end{minipage}
\end{figure}

 We obtain the boundary $\bG$ of the graph $\cG$ by inserting a vertex of valence two at each HR,  the external faces of $\cG$ are incident to these vertices (see an illustration in Figure \ref{fig:boundary}).
 The graph resulting after
the insertion of 2-valent  vertices at each HR is called
pinched ribbon graph \cite{Gurau:2009tz}.

The notion of cut of an edge
as an operation on HERGs has been already introduced. In addition, 
the notions of edge contraction and deletion keep their
ordinary meaning as operations on HERGs. 
(See \cite{rca2} for further details.)

So far, the definition of HERGs is purely combinatorial.
We can certainly comment on the geometrical
meaning of these ribbon structures (this will be
used in section \ref{sect:trans}).

Consider a graph with half-edges (HEG). 
A HEG $\cG$ can be ``cellularly embedded'' in a surface $\Sigma$ with punctures (with boundary circles) in the following sense: 

- Removing all half-edges from $\cG$ we get $\cG'$, a simple
graph, which is then cellularly embedded in $\Sigma$
such that each connected component of $\Sigma\setminus \cG'$ is homeomorphic either to a disc or to discs with holes;

-  Each of the half-edges of $\cG$ is 
embedded in $\Sigma$, respecting of course their incidence relation, and ends on a different puncture 
(but can be on the same boundary circle).  

Then a HERG makes sense as the neighborhood of a HEG cellularly embedded in a punctured surface as defined above.

 We now recall the polynomial invariant introduced in \cite{rca2}.

\begin{definition}[BR polynomial for HERGs \cite{rca2}]
Let $\cG(\cV,\cE,\mf^0)$ be a HERG.  
We define the ribbon graph polynomial of $\cG$ to be
\beq
\cR_{\cG}(X,Y,Z,S,W,T)
=\sum_{A\sset \cG} (X-1)^{\rk(\cG)-\rk(A)}(Y-1)^{n(A)}
Z^{k(A)-F_{\inter}(A)+n(A)} \, S^{C_\partial(A)} \, W^{t(A)}T^{f(A)},
\label{brfla}
\eeq
considered as an element of the quotient of $\mathbb{Z}[X,Y,Z,S,W,T]$ by the ideal generated by $W^2-W$, and 
where $\rk(A)$, $n(A)$, $k(A)$, and $t(A)$ are, respectively, the rank, the nullity, the number of connected components, and the parameter which characterizes the orientability  of $A$ as a surface. If $A$ is orientable, then $t(A)=0$, otherwise, $t(A)=1$. By definition, $\rk(A)=|\cV|-k(A)$ and $n(A)=|\cE(A)|-\rk(A)$. Furthermore, $C_\partial(A)= |\bC(A)|$ is the number of connected components of the boundary of $A$, 
$F_{\inter}(A)= |\cF_{\inter}(A)|$ and $f(A)$ the number of HRs of $A$.
\end{definition}

\begin{theorem}[Contraction and cut on BR polynomial \cite{rca2}]
\label{theo:BRext}
Let $\cG(\cV,\cE,\mf^0)$ be a HERG. Then,
for a regular edge $e$, 
\beq
\cR_{\cG}=\cR_{ \cG\vee e} +\cR_{\cG/e}\,,
\label{retcondel}
\eeq 
for a bridge $e$, we have 
\beq
\cR_{\cG} =(X-1)\cR_{\cG \vee e}+ \cR_{\cG/e}\,;
\label{retbri}
\eeq 
 for a trivial twisted self-loop $e$, the following holds 
\beq
\cR_{\cG}=\cR_{\cG \vee e} + (Y-1)ZW \,\cR_{\cG/e}\,,
\label{retsel2}
\eeq 
whereas for a trivial untwisted self-loop $e$, we have
\beq
\cR_{\cG}=\cR_{\cG \vee e} + (Y-1) \cR_{\cG/e}\,.
\label{retsel1}
\eeq 
\end{theorem}

\vspi 

 The change of variable $S\to Z^{-1}$
 leads to another polynomial $\cR'$ for HERGs 
extending  the BR polynomial. For a graph $\cG$, we write
\beq\label{rpr}
\cR_{\cG}(X,Y,Z,Z^{-1},S,W,T) = \cR'_{\cG}(X,Y,Z,W,T) \,,
\eeq
where $\cR$ is given by  \eqref{brfla}.

\begin{corollary}[Contraction and cut on BR polynomial $\cR'$ \cite{rca2}]
\label{coro:rprim}
Let $\cG(\cV,\cE,\mf^0)$ be a HERG. Then,
for a regular edge $e$, 
\beq
\cR'_{\cG}=\cR'_{ \cG\vee e} +\cR'_{\cG/e}\,,
\qquad 
\cR'_{\cG \vee e} = T^2 \, \cR'_{\cG -e}\,;
\label{retcondelp}
\eeq 
for a bridge $e$, we have 
$\cR'_{\cG/e} = \cR'_{\cG-e} = T^{-2} \, \cR'_{\cG \vee e}$
\beq
\cR'_{\cG} =[(X-1)T^2 +1]  \, \cR'_{\cG/e}\,;
\label{retbrip}
\eeq 
for a trivial twisted self-loop, $\cR'_{\cG-e}=T^{-2} \, \cR'_{\cG \vee e}$ and 
\beq
\cR'_{\cG}=[T^2 + (Y-1)ZW] \,\cR'_{\cG-e}\,,
\label{retselp2}
\eeq 
whereas for a trivial untwisted self-loop, we have
$\cR'_{\cG-e}=T^{-2} \, \cR'_{\cG \vee e}$
and 
\beq
\cR'_{\cG}=[T^2+ (Y-1)] \, \cR'_{\cG-e}\,.
\label{retselp1}
\eeq 
\end{corollary}

\vspi 

Graph operations such as the disjoint union and the one-point-joint ($\cG_1 \sqcup \cG_2$ and $\cG_1 \cdot_{v_1,v_2} \cG_2$, respectively) 
extend  to HERGs \cite{rca}. The product $\cG_1 \cdot_{v_1,v_2} \cG_2$ 
at the vertex resulting from merging $v_1$ and $v_2$ 
keeps its usual sense and respects the cyclic order of all edges and HRs on the previous vertices $v_1$ and $v_2$. The following proposition holds.

\begin{proposition}[Operations on BR polynomials \cite{rca}] 
\label{opBRfla}
Let $\cG_1$ and $\cG_2$ be two disjoint ribbon graphs with HRs, then 
\bea 
\cR_{\cG_1 \sqcup \cG_2}&=&\cR_{\cG_1}
 \cR_{\cG_2},\quad \cR'_{\cG_1 \sqcup \cG_2}= \cR'_{\cG_1}
 \cR'_{\cG_2}\,,
\label{brcups1}\\
\cR'_{\cG_1 \cdot_{v_1,v_2} \cG_2} 
&=&\cR'_{\cG_1} \cR'_{\cG_2} \,,
\label{brcups2}
\eea
for any disjoint vertices $v_{1,2}$ in $\cG_{1,2}$, respectively. 
\end{proposition}

It turns out that $\cR_{\cG_1 \cdot_{v_1,v_2} \cG_2} 
\neq \cR_{\cG_1} \cR_{\cG_2}$. After a closer inspection (see in Appendix \ref{app:semifacto}), we have investigated the breaking terms by giving  more structure to the one-point-joint operation. 
The breaking of this factorization property will a have severe 
consequence on the formulation of a recipe theorem
as we will see in the sequel.

We now introduce a new equivalence relation on HERGs.
   
\begin{figure}[htbp]
\includegraphics[scale=0.5]{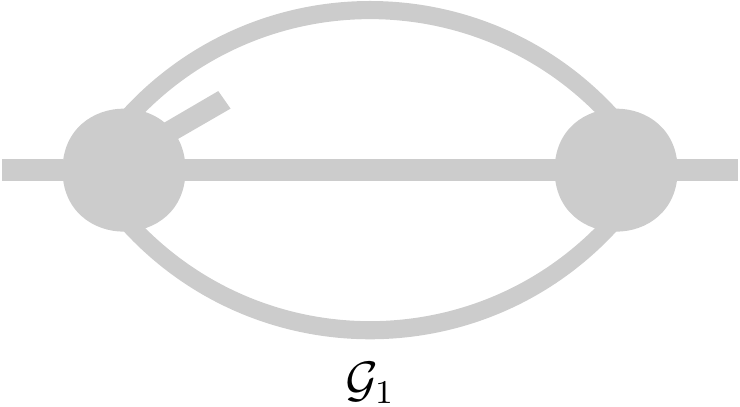}
\hspace{1cm}
  \includegraphics[scale=0.5]{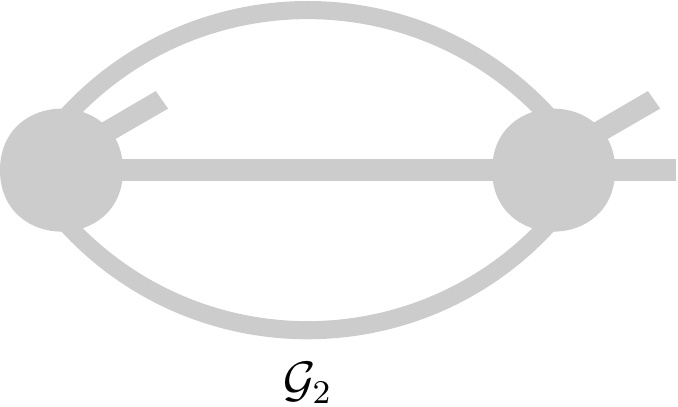}
\hspace{1cm}
 \includegraphics[scale=0.5]{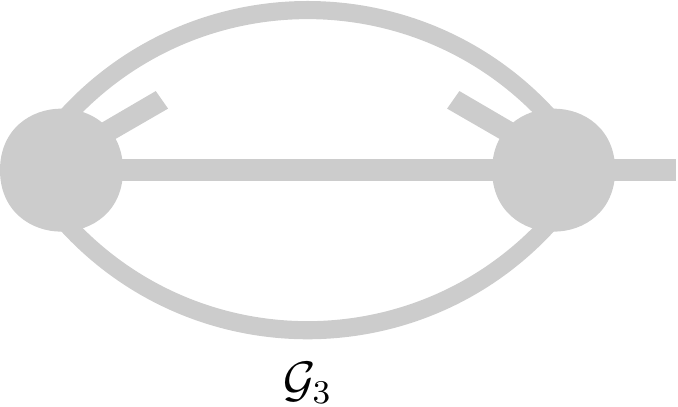}
  \caption{Some HR moves.}
  \label{flag}
\end{figure}

\begin{definition}[HR move operation]
Let $\cG(\cV,\cE,\mf^0)$ be a ribbon graph with HRs. A HR move in $\cG$ consists in removing a HR $f\in\mf^0$ from one-vertex $v$ and placing $f$ either on $v$ or on another vertex such that it is called 

- a HR displacement if the boundary connected component
where $f$ belongs is not modified (see $\cG_1$ and $\cG_2$ in Figure \ref{flag});

- a HR jump if the HR is moved from one boundary connected
component to another one, provided the former remains
a connected boundary component (see $\cG_1$ and $\cG_3$ or $\cG_2$ and $\cG_3$ in Figure \ref{flag}). 
\end{definition}

 HR jumps can modify the boundary graph whereas under HR displacements the boundary graph remains unchanged. 
Nevertheless, the number of connected components of the boundary graph is always preserved under these operations of HR moves.

\begin{definition}[HR-equivalence relation]
We say that two ribbon graphs with HRs $\cG$ and $\cG'$ are HR-equivalent if they are related by a sequence of HR moves. 
\end{definition}

We recall the following statements (Lemma 3 and Lemma 4, 
respectively in \cite{rca2}):

\begin{proposition}
\label{equalmon}
If two ribbon graphs with HRs $\cG$ and $\cG'$ are HR-equivalent, then  for any edge $e$ in $\cG$ and $\cG'$, $\cG\vee e$ and $\cG'\vee e$ are HR-equivalent.
\end{proposition}

\begin{proposition}
\label{equalpol}
 For two HR-equivalent ribbon graphs, $\cG$ and $\cG'$,  $\cR(\cG)=\cR(\cG')$ with $\cR$ the polynomial defined in \eqref{brfla}.
\end{proposition}

It becomes immediate that the definition of the polynomial $\cR$ 
exports on HR-equivalence classes.
 
\begin{definition}[Polynomial for HR-equivalence classes \cite{rca2}]
Let $\cG(\cV,\cE,\mf^0)$ be a ribbon graph with HRs  and $[\cG]$ be its HR-equivalence class. We define the polynomial of $[\cG]$ to be
\bea \label{brequifla}
\cR_{[\cG]}=\cR_{\cG}\,.
\eea
\end{definition}

Thus, for $\cG$ a ribbon graph with HRs, $[\cG]$ its HR-class and $e$ one of its edges, 
$\cR_{[\cG\vee e]}=\cR_{\cG\vee e}$ and 
$\cR_{[\cG/e]}=\cR_{\cG/e}$. We therefore have the next
claim. 

\begin{theorem}[Contraction/cut on BR polynomial on classes \cite{rca2}]
Let $\cG(\cV,\cE,\mf^0)$ be a ribbon graph with HRs and $[\cG]$ be its HR-equivalence class.  Then, for a regular edge $e$, 
\beq
\cR_{[\cG]}=\cR_{ [\cG\vee e]} +\cR_{[\cG/e]}\,,
\label{equiretcondel}
\eeq 
for a bridge $e$, we have 
\beq
\cR_{[\cG]} =(X-1)\cR_{[\cG \vee e]}+ \cR_{[\cG/e]}\,,
\label{equiretbri}
\eeq 
 for a trivial twisted self-loop $e$, the following holds 
\beq
\cR_{[\cG]}=\cR_{[\cG \vee e]} + (Y-1)ZW \,\cR_{[\cG/e]}\,,
\label{equiretsel2}
\eeq 
whereas for a trivial untwisted self-loop $e$, we have
\beq
\cR_{[\cG]}=\cR_{[\cG \vee e]} + (Y-1) \cR_{[\cG/e]}\,.
\label{equiretsel1}
\eeq 
\end{theorem}

To provide a universality property
for $\cR$, let us consider the following expansion of $\cR_{\cG}$:

\bea
&&
\cR_{\cG}(X,Y,Z,S,T,W)
=\sum_{i,j,k,l,m}\cR_{ijklm}(\cG)   (Y-1)^{i} Z^{j} \, S^{k} \, T^{l}\, W^m\label{brfla2},\\
&&
 \cR_{ijklm}(\cG)  :=\sum_{A\sset \cG/ n(A)=i,\,
k(A)-F_{\inter}(A)+n(A)=j,\, C_\partial(A)=k, \,
f(A)=l, t(A)=m} (X-1)^{\rk(\cG)-\rk(A)},
\nonumber
\eea
with $\cR_{ijklm}$  a map from the set $\cG^*$ of isomorphism classes of connected ribbon graphs with HRs to $\mathbb{Z}[X]$
which extends to $\mathfrak{G}$  the set of HR-equivalence classes of isomorphism classes of connected ribbon graphs with HRs.
Then the following statement holds.

\begin{theorem}[Universality of $\cR$ on classes \cite{rca2}] 
\label{theo:univ}
Let $\mathfrak{R}$ be a commutative ring and $x \in \mathfrak{R}$. If a function $\phi: \mathfrak{G} \to \mathfrak{R}$
satisfies 
\bea \label{system}
\phi([\cG]) = \left\{\begin{array}{ll} 
\phi([\cG\vee e]) + \phi([\cG/e]) & {\text{if e is regular}},\\\\
(x-1)\phi([\cG\vee e]) + \phi([\cG/e]) & {\text{if e is a bridge}}.
\end{array} \right. 
\eea
\medskip 
Then there are coefficients $\lambda_{ijklm}\in\mathfrak{R}$, with $i\geq 0$, $0\leq k\leq i+1$, $l\geq 0$, $0\leq m\leq 1$ and $0\leq j \leq i+1$  such that 
\beq \label{univsum}
\phi([\cG]) = \sum_{i,j,k,l,m} \lambda_{ijklm} \cR_{ijklm}(x).
\eeq
\end{theorem}

\vspi 
\noindent{\bf Chord diagrams.}
So far, all proofs of universality theorems for topological 
polynomials on ribbon graph rest on the ``projection'' of
one-vertex ribbon graphs onto the so-called chord-diagrams
\cite{bollo,joan,rca2}.  
These diagrams will  play a crucial role in the 
proof of the following recipe theorems as well. In 
a way useful to the present context, such a notion of
chord diagrams and its main properties must be  recalled.

\begin{definition}[Chord diagrams \cite{rca2}]
$\bullet$ A half-chord on a chord diagram is a segment attached to a unique point on its circle. 

$\bullet$ An (open) chord diagram is a chord diagram in sense of \cite{bollo} with a (nonempty) set of half-chords. In the case where this set is empty, it becomes BR chord diagram.

$\bullet$ A signed (open) chord diagram is an (open) chord diagram with an assignment of a sign ``\tw'' or not to each chord.
\end{definition}

\begin{figure}[htbp]
\includegraphics[scale=1]{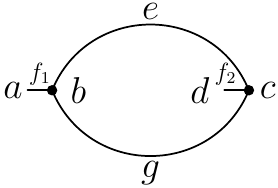}
\caption{Two-vertex ribbon graph with HRs.}
\label{twovert}
\end{figure}

\begin{figure}[htbp]
   \includegraphics[scale=0.8]{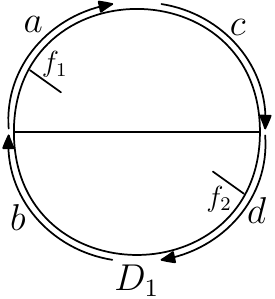}
\quad 
 \includegraphics[scale=0.8]{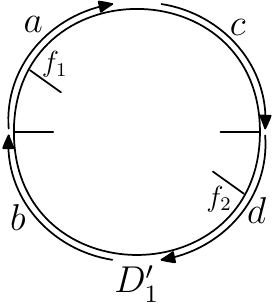}
\quad
 \includegraphics[scale=0.8]{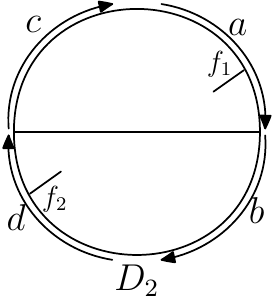}
\quad
 \includegraphics[scale=0.8]{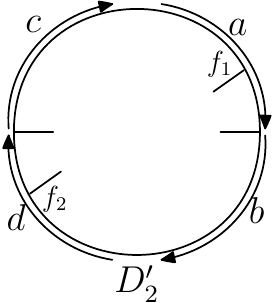}
  \caption{Related chords diagrams $D_1$, $D'_1$, $D_2$, $D'_2$.}
  \label{oprotate}
\end{figure}

\begin{figure}[htbp] 
   \includegraphics[scale=0.8]{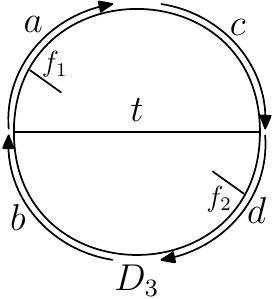}\quad
 \includegraphics[scale=0.8]{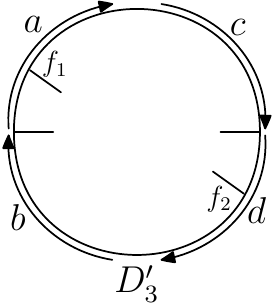}\quad
\includegraphics[scale=0.8]{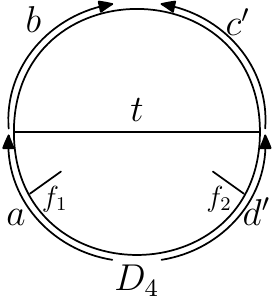}\quad
 \includegraphics[scale=0.8]{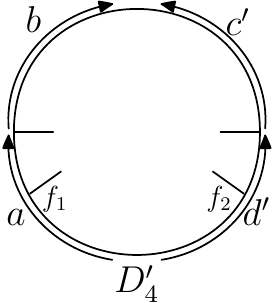}
  \caption{Related chords diagrams $D_3$, $D'_3$, $D_4$, $D'_4$.}
  \label{optwiste}
\end{figure}

Given a one-vertex HERG, we assign  twisted ribbon edges 
to  chords coined by ``$\tw$'' in the corresponding chord
diagram.  

Consider a two-vertex (embedded) half-edged graph as given in Figure \ref{twovert}. 
 Contracting either $e$ or $g$ yields a one-vertex graph 
which can be mapped onto a chord diagram. 
Note that $e$ or $g$ might correspond to twisted edges. 
The different possibilities are summarized by Figures \ref{oprotate}
($D_1$ and $D_2$ if both $e$ and $g$ are not twisted) and \ref{optwiste} ($D_3$ and $D_4$ if one of $e$ or $g$ are twisted). 
Then from each of these configurations, one cuts the chord
coming from the edge about which we rotate or twist (this gives the configurations
denoted by $D'_{1,2,3,4}$). Detailed  explanations can be found
in the companion paper \cite{rca2}.

Two signed (open) chord diagrams are said to be related by a rotation about the chord $e$ if they are related as $D_1$ and $D_2$ in Figure \ref{oprotate}, and related by a twist about $e$, if they are related as $D_3$ and $D_4$ in Figure \ref{optwiste}.  We now give the definitions of $R$-equivalent diagrams and the sum of two chord diagrams.

\begin{definition}[$R$-equivalence relation \cite{bollo}]
Two diagrams or signed  diagrams $D_1$ and $D_2$ are $R$-equivalent if and only if they are related by a sequence of rotations and twists. We write  $D_1\sim D_2$.
\end{definition}
\begin{definition}[Sum of  diagrams \cite{bollo}]
The sum of two  diagrams or signed  diagrams $D_1$ and $D_2$ is obtained by choosing a point $p_i$ (not the end-point of a chord or a half-chord) on the boundary of each $D_i$, joining the boundary circles at these points and then deforming the result until it is again a circle.
\end{definition}

The sum of diagrams can be performed in several ways, if  the $p_i$ are chosen differently. However, all of them are $R$-equivalent. 
The proof of the following statement can be found in \cite{rca2}
(Lemma 1) which is an adjustment of a similar lemma in \cite{bollo}.

\begin{lemma}\label{lem3}
If two diagrams $D$ and $D'$ are both sums of diagrams $D_1$ and $D_2$, then they are $R$-equivalent.
\end{lemma}

The next important notion is that of canonical chord diagrams. Given $i\geq 0$,  $0\leq 2j\leq i$, $0\leq k\leq i+1$, $l\geq0$ and $0\leq m\leq 2$, let $D_{i,j,k,(s;l_1,\cdots ,l_q),m}$ be the chord diagram consisting of $i$ chords, $j$ pairs of chords intersecting each other, $k$ connected components of the boundary of this diagram, $l$ half-chords ($l = s + \sum_{p=1}^q l_p$) specifically arranged and $m$ negative chords (or twisted chords) intersecting no other chords ($i-2j-m$ is the number of positive chords intersecting no other chords). This diagram is drawn such that a number $l-s$ of half-chords is partitioned in $(l_p)_{p=1,\cdots, q}$, positive chords intersecting no other chords (we call these isolated chords) and $s$ is the rest of the half-chords. All these chords and half-chords are arranged on the circle diagram (see an illustration for $D_{4,1,2,(3;1),1}$ and  $D_{5,1,2,(0;1,2),1}$):

\begin{figure}[htbp]
 \includegraphics[scale=0.8]{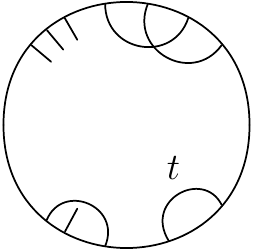}
\hspace{1cm}
 \includegraphics[scale=0.8]{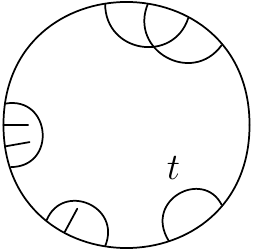}
\caption{Canonical diagrams: $D_{4,1,2,(3;1),1}$ and $D_{5,1,2,(0;1,2),1}$.}
  \label{canoni}
\end{figure}

It becomes straightforward to perform a quotient on the space of
chord diagrams with respect to the HR-equivalence relation.
One gets HR-equivalent (canonical) chord diagrams. 
At this stage, the partition of $l$ half-chords becomes
superfluous to emphasize. One can dispose the $l$ flags on the chord
diagram as desired and must only respect the number of connected components of the boundary. It can proved that Lemma 
\ref{lem3} still holds in that context.

\section{Recipe theorems}
\label{sect:Recipe}

This section discusses recipe theorems for the BR polynomial  on classes of HR-equivalent ribbon graphs.  
Then, only in this section,  we call 
HERG or ribbon graph with half-edges (or with HR), a class of HR-equivalent class. The same remark holds for HR-classes 
of (canonical) chord diagrams which are simply refereed  
to (canonical) chord diagrams. 
For simplicity,  
we drop all brackets hence, from now on,
$[\cG]$ will be denoted simply $\cG$, 
and $D_{ijklm}$ stands for the HR-equivalent class 
canonical diagram $[D_{i,j,k,[l],m}]$. 

Any (open) diagram is related to a canonical signed chord diagram $D_{ijklm}$ ($0\leq m\leq 2$) that consists of $i-2j-k-m$ positive disjoint chords $j$ pairs of intersecting positive chords, $k$  connected components of the boundary graph, $m$ negative disjoint chords  and $l$ half-chords.
In this case, we observe that: 
\bea
\mathcal{R}(D_{10000};x,y,z,s,t,w) = (y-1) + zst^2,\nonumber\\
\mathcal{R}(D_{10001};x,y,z,s,t,w) = (y-1)zw + zst^2,\nonumber\\
\mathcal{R}(D_{21000};x,y,z,s,t,w) = (y-1)^2z^2 + 2(y-1)(zst)^2 + zst^4,\nonumber\\
\mathcal{R}(D_{001l0};x,y,z,s,t,w) = zst^l. \label{poly:canonical}
\eea

The polynomial $\mathcal{R}$ is not multiplicative with respect to the one-point-joint
 of ribbon graphs with half-edges (see Appendix \ref{app:semifacto}
for a complete development on this issue). In order to find a recipe, this factorization property must be replaced by other useful relations.
We have identified  weaker conditions that $\cR$ satisfies
which will ensure the existence of a recipe in these cases.  

 In the following, we denote by $M(A)$ the monomial associated with 
 $A\sset \cG$ which we recall
\bea
M(A)(x,y,z,s,t,w)=(x-1)^{\rk(\cG)-\rk(A)}(y-1)^{n(A)}
z^{k(A)-F_{\inter}(A)+n(A)} \, s^{C_\partial(A)} \, w^{t(A)}t^{f(A)}.
\eea
Let $\gamma$ be the function defined on the set $\{0,1,2\}$ by: 
\bea \label{eq:gamma}
\left\{\begin{array}{ll} 
\gamma(0)=0, \\
\gamma(1)=\gamma(2)=1.
\end{array} \right. 
\eea
 \begin{lemma}\label{lem:decM}
Consider a canonical chord diagram with half-chords  $D_{ijklm}.$ Then, the following decomposition is valid:
\bea 
M(D_{ijklm})(x,y,z,s,t,w)=(M(D_{10000}))^{i-2j-m}(M(D_{001l0}))^k(M(D_{10001}))^{m}\times \nonumber \\(M(D_{21000}))^{j}(x,y,z,s,t^{\frac{1}{k}},w) \,.\label{onepointjoint}
\eea
\end{lemma} 
\proof  We have $M(D_{ijklm})(x,y,z,s,t,w)=(y-1)^iz^{1-F_{\inter}(D_{ijklm})+i}s^kt^lw^{\gamma(m)}$ with $\gamma$ defined as in  \eqref{eq:gamma} and $F_{\inter}(D_{ijklm})$ $=i-2j-m-k+1.$ Then, $M(D_{ijklm})=(y-1)^iz^{2j+m+k}s^kt^lw^{\gamma(m)}$. Furthermore,
evaluating  $M(D_{10000})(x,y,z,s,t^{\frac{1}{k}},w)=y-1$,
 $M(D_{001l0})(x,y,z,s,t^{\frac{1}{k}},w)$ $=zst^{\frac{l}{k}}$, 
$M(D_{10001})(x,y,z,s,t^{\frac{1}{k}},w)$ $=(y-1)zw$ and $M(D_{21000})$ $(x,y,z,s,t^{\frac{1}{k}},w)=(y-1)^2z^2$. This ends the proof of \eqref{onepointjoint}.

\qed

From the definition of $\mathcal{R}$, 
\bea
\mathcal{R}(D_{ijklm})(x,y,z,s,t,w)=\sum_{D_{i'j'k'l'm'}\sset D_{ijklm}} M(D_{i'j'k'l'm'})(x,y,z,s,t,w), \label{nonepointjoint}
\eea
because any subgraph of $D_{ijklm}$ is of the form $D_{i'j'k'l'm'}$ with $i'\leq i$.
We apply \eqref{nonepointjoint} in order to find the following relations:
\bea
\mathcal{R}(D_{10000})&=&M(D_{10000}) + M(D_{00120}),\cr
\mathcal{R}(D_{21000})&=&M(D_{21000}) + 2M(D_{10000})(M(D_{00110}))^2 + M(D_{00140}),\cr
\mathcal{R}(D_{10001})&=&M(D_{10001}) + M(D_{00120}),\cr
\mathcal{R}(D_{001l0})&=&M(D_{001l0}).
\label{diagram}
\eea
These relations can be inverted to find
\bea
M(D_{10000})&=& \mathcal{R}(D_{10000}) - \mathcal{R} (D_{00120}),\cr
M(D_{10001})&=& \mathcal{R}(D_{10001}) - \mathcal{R}(D_{00120}),\cr
M(D_{21000})&=& \mathcal{R}(D_{21000}) +2(\mathcal{R}(D_{10000}) - \mathcal{R}(D_{00120}))(\mathcal{R}(D_{00110}))^2 + \mathcal{R}(D_{00140}),\cr
M(D_{001l0})&=&\mathcal{R}(D_{001l0}).
\label{dia}
\eea
We immediately discover  that,  from Lemma \ref{lem:decM},
the existence of such relations \eqref{diagram} allows us to 
write in return  a generic monomial $M(D_{ijklm})$ in terms of $\cR (D_{i'j'k'l'm'})$,
for $i',j',k',l'$ and $m'$ well chosen. 
This means that the system \eqref{nonepointjoint} can be inverted. 
In general $M(D_{ijklm})$ is a non linear function of $\cR(D_{ijklm})$. 
For the simplest chord diagrams \eqref{dia}, we observe that 
the relation is still linear.   

	Replacing the relations \eqref{diagram} in \eqref{onepointjoint} and in \eqref{nonepointjoint}, we obtain $\mathcal{R}(D_{ijklm})$ as function of 
 $\mathcal{R}(D_{10000})$,  $\mathcal{R}(D_{10001})$,  $\mathcal{R}(D_{21000})$ and  $\mathcal{R}(D_{001l0})$.

 Let $\mathcal{M}$ be a minor closed subset of HERGs containing all HERGs on two vertices (this notion extends \cite{iain} to HERGs). Let $F$ be a map from $\mathcal{M}$ to a commutative ring $\mathfrak{R}$ with unity.
Suppose that $F$ satisfies the relation:
\bea
&&
F(D_{ijklm})=\sum_{D_{i'j'k'l'm'}\sset D_{ijklm}} \alpha^{i'-j'+m'+\frac{l'}{k'}}M(D_{i'j'k'l'm'}) \label{nonepointjoint2}\\
&&=\sum_{i'\leq i} (M(D_{10000}))^{i'-2j'-m'}(M(D_{001l'0}))^{k'}(M(D_{10001}))^{m'}(M(D_{21000}))^{j'}(x,y,z,s,t^{\frac{1}{k'}},w),\nonumber
\eea
with the conditions: 
\bea\label{bapteme1}
M(D_{10000})&=&\alpha^{-1} [F(D_{10000}) -F (D_{00120})],
\cr
M(D_{10001})&=&\alpha^{-1} [F(D_{10001}) - F(D_{00120})],
\cr
M(D_{21000})&=&\alpha^{-1} [F(D_{21000}) + (F(D_{10000}) - F(D_{00120}))(F(D_{00110}))^2 + F(D_{00140})],
\cr
M(D_{001l0})&=&\alpha^{-1}F(D_{001l0}).
\eea
From  \eqref{bapteme1}, we are guaranteed that
once again the monomials $M(D_{ijklm})$ can be written in terms
of $F(D_{i'j'k'l'm'})$, for $i',j',k',l'$ and $m'$ well chosen. 
Under conditions \eqref{nonepointjoint2} and \eqref{bapteme1}, the recipe 
theorem implies that $F$ is strongly related to $\cR$.

\begin{theorem}[First recipe theorem for $\cR$]
\label{theo:recipe}
Let $\mathcal{M}$ be a minor closed subset of HERGs containing all HERGs on two vertices. Let $F$ be a map from $\mathcal{M}$ to a commutative ring $\mathcal{R}$ with unity. Let $s=F(D_{21000})$, $q=F(D_{10000})$, $r=F(D_{10001})$ and $s_l=F(D_{001l0})$ and suppose there are elements $\alpha$, $x$, $u$, $v$, $w$, $o$ $\in \mathfrak{R}$ with $\alpha$ a unit such that:
\begin{enumerate}
\item \bea \label{system1}
F(G) = \left\{\begin{array}{ll} 
F(G \vee e) + F(G/e) & {\text{if e is  regular}},\\\\
(x-1) \, F(G \vee e) + F(G/e) & {\text{if e is a bridge}}.
\end{array} \right. 
\eea
\item  $F(G\sqcup H)=F(G)F(H)$ where $G$ and $H$ are embedded bouquets with half-edges
and  $F$ satisfies also \eqref{nonepointjoint2} and \eqref{bapteme1};
\item $F(E)=\alpha^n$ if $E$ is an edgeless graph  with $n$ vertices
without HRs;
\item  $(q- s_2)^2u^2=\alpha[s-2\alpha^{-2}(q-s_2)s_1^2 -  s_4]$, and $(q-s_2)uw=r- s_2$, and $s_l=\alpha(uvo^l)$, and also $w=w^2$. 
Then
\bea
F(G)=\alpha^{k(G)}\mathcal{R}(G;x,\alpha^{-1}(q-s_2+\alpha),u,v,w,o),\label{recipe}
\eea 
where $k(G)$ is the number of components of $G$.
\end{enumerate}
\end{theorem}
\proof The proof of the recipe theorem follows the lines of the universality proof, because as previously
noticed, it is a very similar statement. Thus, we will 
proceed by two inductions, on the number of chords
in chord diagrams representing an embedded bouquet graph
and then on the number of non-loop edges in a general graph. 

From item (1), we obtain
\bea
F(D_1) - F(D'_1) = F(D_2) - F(D'_2), \label{related1} \mbox{ and }\\
F(D_3) - F(D'_3) = F(D_4) - F(D'_4), \label{related2}
\eea
where  $D_1$, $D_2$, $D_3$ and $D_4$ are related as in Figures \ref{oprotate} and \ref{optwiste} and $D'_i=D_i\vee e$.

 Since the polynomial $\mathcal{R}$ satisfies  \eqref{related1} and \eqref{related2}, $F'=F - \alpha \mathcal{R}$ also satisfies the same equations. Using items (2), (3) and (4) of Theorem \ref{theo:recipe}, the relation \eqref{recipe} holds in the case of chord diagrams with 0 chord. 

	Assume by induction that $F(D) = \alpha \mathcal{R}(D)$ for any signed chord diagram with fewer than $n$ chords. In this case $F'$ vanishes on signed chord diagrams with fewer than $n$ chords.
	Using  \eqref{related1} and \eqref{related2}, we have $F'(D)=F'(D_{njklm})$ where $D$ is related to a canonical diagram $D_{njklm}$. Now from \eqref{nonepointjoint2} and \eqref{bapteme1} in  item (2), we deduce:
\bea
F(D_{njklm})=\alpha\mathcal{R}(D_{njklm};x,\alpha^{-1}(q-s_2+\alpha),u,v,w,o)
\eea
and
\bea
F'(D)=F'(D_{njklm})=0.
\eea
	By induction, $F'(D)=0$ on any signed chord diagram $D$. Finally the result holds on any rosette ribbon graph. 
	From  item (1), the result becomes true on any ribbon graph.

\qed

The fact that $\mathcal{M}$ contains all ribbon graphs on two vertices is not important in this proof. The important thing is to require that  $F$ satisfies \eqref{related1} and \eqref{related2} for a chain of rosettes ending in a canonical chord diagram. This gives the so-called
``low fat'' recipe theorem \cite{joan}. Thus, 
we have also a low fat recipe for $\cR$. First, we must introduce
another minor set. A minor closed set $\mathcal{M}$ of HERGs is said closed  under chord operations whenever $D\in\mathcal{M}$ and $D\sim D_{ijklm},$ then there is a finite sequence $D=D_1...D_n=D_{ijklm}$ with $D_i\in\mathcal{M}$ and $D\sim D_{i+1}$ for all $i$.

\begin{theorem}[The ``low fat'' recipe theorem for $\cR$]
\label{theo:lowrecipe} 
Theorem \ref{theo:recipe} holds with ``Let $\mathcal{M}$ be a minor closed subset of HERGs containing all HERGs on two vertices and let $F$ map $\mathcal{M}$ to a commutative ring $\mathfrak{R}$ with unity'', replaced by ``Let $\mathcal{M}$ be a minor closed subset of HERGs closed under chord operations that contains $D_{10000}$, $D_{10001}$, $D_{21000}$, $D_{001l0}$ and let $F$ map $\mathcal{M}$ to a commutative ring $\mathfrak{R}$ with unity, such that $F$ satisfies 
 \eqref{related1} and \eqref{related2} whenever the $D_i's$ are related as in Figures \ref{oprotate} and \ref{optwiste}.''
\end{theorem}
We also have:
\begin{corollary}
If $F$, $\mathcal{M}$, $\mathfrak{R}$ satisfy conditions of Theorem \ref{theo:lowrecipe}, with both $q- s_2$ and $r- s_2$ being unit of $\mathfrak{R}$, then $w=1$, and thus $F$ does not discern orientation by the presence or absence of a single idempotent element.
\end{corollary}
\proof
From the second equation in item (4), we obtain $\frac{r- s_2}{q- s_2}=uw=uw^2=w\frac{r- s_2}{q- s_2}$.  Since $q- s_2$ and $r- s_2$ are units then we infer $w=1$. Then $F(G)=\alpha^{k(G)}\mathcal{R}(G;x,\alpha^{-1}(q-s_2+\alpha),u,v,1,o)$.

\qed

\bigskip

There exists another way to introduce a recipe that we now detail. 
Instead of imposing  that $F$ satisfies \eqref{nonepointjoint2} and \eqref{bapteme1}, we can impose another type of condition. 
 This will give rise to another recipe. To be precise the way to evaluate
the function is the same, but the full set conditions that 
must satisfy $F$ is now different. In a very interesting way, 
the following analysis can be also applied to the BR polynomial $R$
on ribbon graphs (without HRs) and we can recover a recipe theorem  without mentioning a factorization property of the function $F$
with respect to the one-point-joint operation. 
We first need a series of preliminary results and will introduce a specific terminology now.

\begin{figure}[h]
\centering
\begin{minipage}[t]{.8\textwidth}
\centering
\includegraphics[angle=0, width=1.6cm, height=1.5cm]{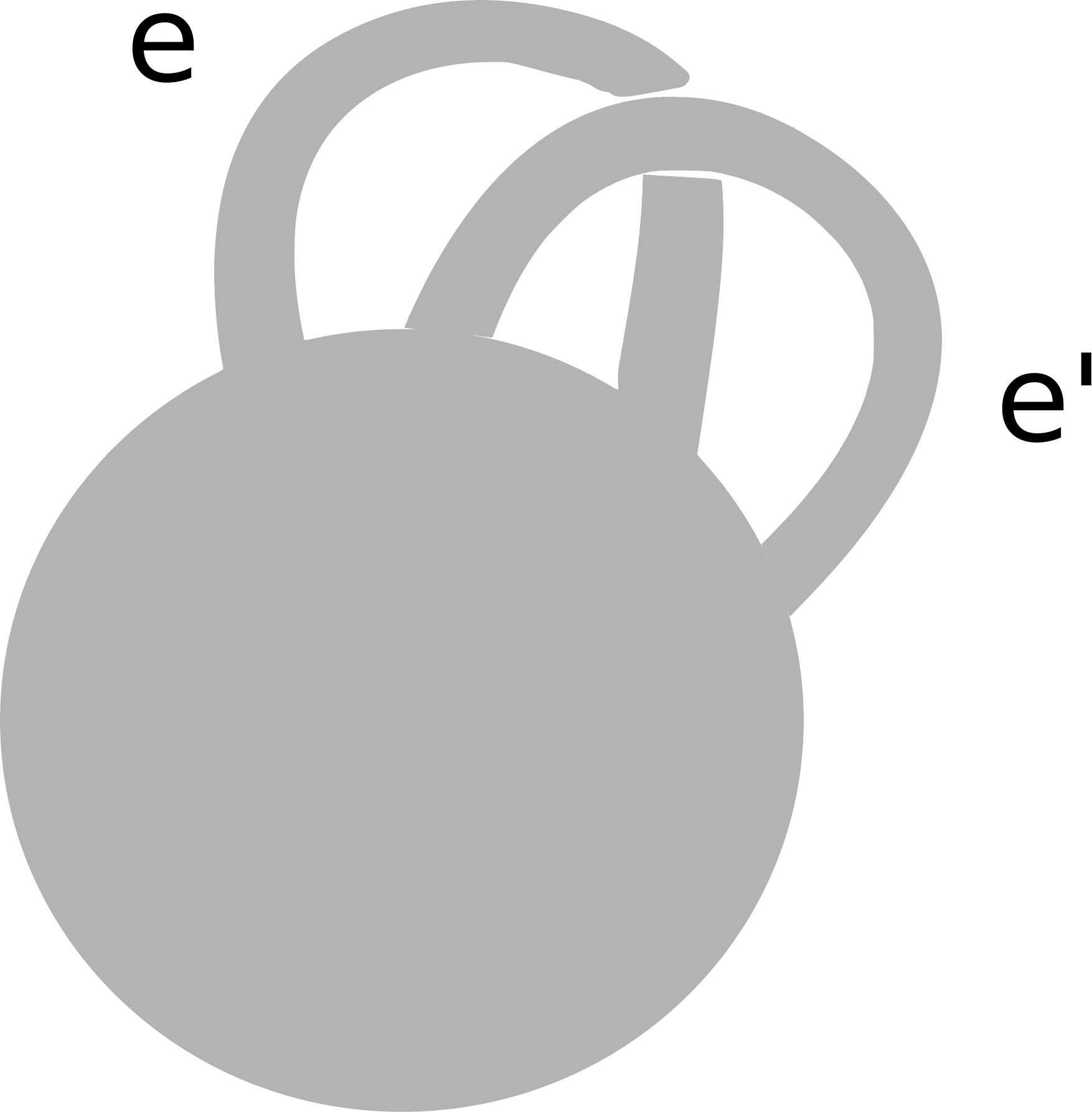}
\caption{ {\small A genus loop.
}}
\label{fig:genusloop}
\end{minipage}
\end{figure}

 A ``genus loop'' of a HR-equivalent ribbon graph $\cG$ is defined by exactly two untwisted loops $e$ 
and $e'$ crossing each other on a vertex and which do not 
cross any other loops at that vertex (see Figure \ref{fig:genusloop}).

\begin{lemma}[Genus loop evaluation] \label{lem:cross}
Let $\cG$ be an one-vertex ribbon graph having a genus
loop defined by $e$ and $e'$. Then, the following recurrence relation is obeyed
\bea
\cR(\cG) =  \cR(\cG \vee e) + ((Y-1)Z)^2 \cR(\cG/e/e') 
+ (Y-1) \cR( (\cG/e) \vee e') \,.
\eea
\end{lemma}
\proof 
The proof is quite reminiscent of the ordinary contraction/cut rule
with however   modifications that we will emphasize. 
From the ordinary partition of the set of spanning subgraphs of $\cG$
between cutting subgraphs containing $e$ and those
which do not,  it is clear that the sum of monomials of subgraphs which do not contain $e$ are directly mapped to $\cR(\cG \vee e)$. 
Now we focus on the rest of subgraphs which contain $e$ and further partition
them into those containing $e'$ and those which do not. 

A subgraph $A$, such that $e,e' \in A$, must be mapped
onto $A/e/e' \sset \cG/e/e'$, by  contracting successively
$e$ and $e'$. To compare the monomials $M(A)$ and
$M(A/e/e')$, we need the following relations (all the rest
of the combinatoric numbers are constants)
\bea
n(A/e/e') = n(A/e) = n(A) -2 .
\eea 
This gives a relation between the two monomials
as $M(A) = ((Y-1)Z)^2 M(A/e/e')$. 

On the other hand, to a subgraph $A$ such that $e\in A$ and $e'\notin A$, we assign $(A/e) \vee e'$. The monomials $M(A)$ 
and $M((A/e) \vee e')$ can be compared as well. We have the basic
relations:
\bea
&&
n((A/e)\vee e') = n(A/e) = n(A) -1 , \cr
&&
r((\cG/e) \vee e') - r((A/e)\vee e') = 0= 
r(\cG/e) - r(A/e) -1 , \cr
&&
r(\cG/e) - r(A/e) = r(\cG) - r(A) + 1, \cr
&&
k((A/e)\vee e') = k(A/e) = k(A) +1. 
\eea 
Therefore, $M(A) = (X-1)^{-1} (Y-1) M(A/e) 
= (X-1)^{-1} (Y-1) (X-1) M((A/e)\vee e') 
 = (Y-1) M((A/e)\vee e')$. 
Summing over all contributions gives 
the result. 

\qed

The following statement holds. 

\begin{theorem}[Second recipe theorem for $\cR$]
\label{theo:recipe2}
Let $\mathcal{M}$ be a minor closed subset of 
HERGs containing all HERGs on two vertices. Let $F$ be a map from $\mathcal{M}$ to a commutative ring $\mathfrak{R}$ with unity. Let $s=F(D_{21000})$, $q=F(D_{10000})$, $r=F(D_{10001})$ and $s_l=F(D_{001l0})$ and suppose there are elements $\alpha$, $x$, $y$, $z$, $u$, $v$, $w$, $o$ $\in \mathfrak{R}$ with $\alpha$ a unit such that:
\begin{enumerate}
\item \bea \label{systg}
F(G) = 
\left\{\begin{array}{ll} 
F(G \vee e) + F(G/e) & {\text{if e is  regular}},\\
(x-1) \, F(G \vee e) + F(G/e) & {\text{if e is a bridge}},\\
F(G \vee e) + \alpha^{-1}(y-1)F(G/e) & {\text{if e is  a trivial untwisted loop,}}\\
F(G \vee e) + (y-1)zw F(G/e) & {\text{if e is a trivial twisted loop}}, \\
 F(\cG \vee e) + ((y-1)z)^2 F(\cG/e/e') \\
\qquad \quad\; +\alpha^{-1} (y-1) F( (\cG/e) \vee e')
& {\text{if $(e,e')$ defines a genus loop}}, 
\end{array} \right. 
\eea
\item  $F(G\sqcup H)=F(G)F(H)$ where $G$ and $H$ are embedded bouquets with half-edges.
\item $F(E)=\alpha^n$ if $E$ is an edgeless graph  with $n$ vertices
without HRs;
\item  $(q- s_2)^2u^2=\alpha[s-2\alpha^{-2}(q-s_2)s_1^2 -  s_4]$, and $(q-s_2)uw=r- s_2$, and $s_l=\alpha(uvo^l)$, and also $w=w^2$. 
Then
 \bea
F(G)=\alpha^{k(G)}\mathcal{R}(G;x,\alpha^{-1}(q-s_2+\alpha),u,v,w,o),\label{recip2}
\eea 
where $k(G)$ is the number of components of $G$.
\end{enumerate}
\end{theorem}

Before proving this theorem we shall need a particular 
relation satisfied by $F$ on canonical chord diagrams. Here we will denote a  canonical chord diagram $D_{i,j,k,l,m}$ because we will perform operations on the indices.

\begin{lemma}\label{lemma:fequaldiag}
For all $i,j,k,l$ and $m$ for an arbitrary chord diagram $D_{i,j,k,l,m}$,
and let $F$ a function on HERGs satisfying the condition (1)-(3) and $F(D_{0,0,1,l,0}) = \alpha \cR(D_{0,0,1,l,0})$
then we have
\bea \label{equ:fequaldiag}
F(D_{i,j,k,l,m}) = \alpha \cR(D_{i,j,k,l,m}).
\eea
\end{lemma}

\proof  In order to prove our claim, we will use an algorithm
which will reduce the number of chords and the 
complexity of the chord diagrams. From simpler cases,
we will be able to conclude.

If $m>0$, there exists a negative chord $e$ and we have
\bea\label{fm}
F(D_{i,j,k,l,m} ) &=& F( D_{i,j,k,l,m} \vee e) + (y-1)zw F(D_{i,j,k,l,m}/e) 
\crcr
&=& F( D_{i-1,j,k + \varepsilon,l+2,m-1} ) + (y-1)zw F(D_{i-1,j,k,l,m-1}) 
\eea
where $\varepsilon=0,1$ depends on the type of chord diagram. 
The property \eqref{fm} reduces the number of negative
chords up to 0 (we can use it twice if $m=2$). We note also
that the total number of chords $i$  decreases. After the
procedure, one gets a sum of (2 or 4) terms involving new
$D_{i',j',k',l',0}$. Note that a similar relation holds exactly
for $\cR$. 
 We can therefore concentrate on the evaluation of such 
canonical diagram  $D_{i,j,k,l,0}$.

(A) For chords which do not intercept any other chords 
and which encloses $\ell$ HRs, we use the  relation in \eqref{systg}
corresponding to trivial untwisted loops. The number of these chords is $i-2j$.
We will contract/cut these chords until none will
be left. Considering that the loop 
$e$ has $\ell$ (possibly 0) HRs, it can generate $\epsilon=0$ or $1$ connected component of the boundary. 
For such a chord $e$ enclosing $\ell\geq 0$ chords, we write
\bea
F(D_{i,j,k,l,0}) = F(D_{i,j,k,l,0}\vee e)
 + \alpha^{-1} (y-1) F(D_{i,j,k,l,0}/e). 
\eea
Note that the first term  $F(D_{i,j,k,l,0}\vee e) =F(D_{i-1,j,k+\varepsilon,l+2,0})$, $\varepsilon=\pm 1,0$, does not cause any trouble:
the number of chords decreases and we will perform again (A). 
The possible terms generated are of the same
form as for the initial $D_{i,j,k,l,0}$ and the following discussion can 
be applied for $F(D_{i-1,j,k+\varepsilon,l+2,0})$. 
 We then concentrate on the second term. Using in particular (3),
and the function $\epsilon(\ell)=0$, if $\ell =0$, and
$\epsilon(\ell)=1$, if $\ell >0$,   
we write this term as 
\bea
&&
\alpha^{-1}(y-1) F(D_{i,j,k,l,0}/e) = 
\alpha^{-1}(y-1) F(D_{i-1,j,k-\epsilon(\ell),l-\ell,0})F(D_{0,0,\epsilon(\ell),\ell,0})
\crcr
&&
 = 
(y-1) F(D_{i-1,j,k-\epsilon(\ell),l-\ell,0})\cR(D_{0,0,\epsilon(\ell),\ell,0})\,.
\eea
The last simplification occurs by definition $F(D_{0,0,1,l,0})
 = \alpha \cR(D_{0,0,1,l,0})$, $\forall l\geq 0$. 
Thus the evaluation of $F$ on $D_{i,j,k,l,0}$ with chords enclosing
possibly  HRs revolves into a sum of terms involving the evaluation 
of $F$ on $D_{i',j',k',l',0}$ with fewer number of chords
times some power of $(y-1)$ and then products of
 terms involving $\cR$. 
Evaluating $\cR$ itself on  the same $D_{i,j,k,l,0}$ 
will give a similar expression apart from the
evaluation of $\cR$ on $D_{i',j',k',l',0}$. Then apply again (A).

(B) The last type of chords are genus loops. We adopt the same
strategy as above. The argument is lengthier but we can 
also prove that, using the particular relation adapted to
these genus loops in \eqref{systg} (or Lemma \ref{lem:cross}), we can recast the 
evaluation of $F(D_{2j,j,k,l,0})$ a sum of terms of 
simpler diagrams. We will simply give here the list
of arguments which will allow one to achieve the proof. 

- Consider $(e,e')$ a genus loop. Use the double recurrence
relation and find that the evaluation of $F$ on $D_{2j,j,k,l,0}$
involves the evaluation of $F$ onto

(a) $D_{2j,j,k,l,0}\vee e$ is of the form $D_{2j-1,j-1,k+\epsilon,l+2,0}$ possesses a untwisted loop. One must reduce this untwisted loop
using  (A) and find a sum of 
terms where $F$ evaluates as a sum of simpler terms 
on a unique canonical diagram $D_{i',j',k',l',0}$. 

(b) $D_{2j,j,k,l,0}/e/e'$ which is of the form $D_{2j-2,j-1, k, l,0}$
and we go back to the procedure (B). 

(c) $(D_{2j,j,k,l,0}/e) \vee e'$ which is of the form 
$D_{2j-2,j-1,k+\epsilon, l+1,0} \sqcup D_{0,0,1,1,0}$,
then we use the fact that $\alpha^{-1} F(D_{0,0,1,1,0})=\cR(D_{0,0,1,1,0})$,
to recast the evaluation of $F$ on this sector
to a single chord diagram $D_{2j-2,j-1,k+\epsilon, l+1,0}$, 
and go back to (B).

At the end of the algorithm, $F(D_{i,j,k,l,m})$ is a sum of
terms which involves products of the form $\alpha^{-1} f(y-1,z,w) F(D_{0,0,\epsilon(l'),l',0}) \prod_{\kappa \in K}\cR(D_\kappa)$, where
 $l' \geq 0$,  $f(y-1,z,w)$ is a polynomial function which is a product of different contributions of special edges which 
have been contracted or cut and $K$ is a family of 
chord diagrams of the form $D_{\kappa}= D_{0,0,\epsilon(l),l,0}$. 
If we apply the same procedure to $\cR(D_{i,j,k,l,m})$
each term in this expansion will be of the form
 $f(y-1,z,w) \cR(D_{0,0,\epsilon(l'),l',0})\prod_{\kappa \in K}\cR(D_\kappa)$. We conclude that $F(D_{i,j,k,l,m})=\alpha
\cR(D_{i,j,k,l,m})$, since one has $F(D_{0,0,\epsilon(l),l,0}) = \alpha \cR(D_{0,0,\epsilon(l),l,0})$,
for all $l\geq 0$. 

\qed

The proof of our theorem can be now completed. 

\proof[Proof of the Theorem \ref{theo:recipe2}]
Using Lemma \ref{lemma:fequaldiag}, the relation \eqref{recip2} holds in the case of canonical chord diagrams with half-chords. 

 Using the two first equations in item (1), $F$ satisfies  \eqref{related1} and \eqref{related2} for any chord diagrams $D_1$, $D_2$, $D_3$ and $D_4$ related as in Figures \ref{oprotate}
and \ref{optwiste} and $D'_i=D_i\vee e$. Hence $F'=F - \alpha \mathcal{R}$ also satisfies the same equations since $\mathcal{R}$ satisfies them.

We prove that this theorem holds on any chord diagram using step by step the proof in Theorem \ref{theo:recipe} and replacing the relations \eqref{nonepointjoint2} and \eqref{bapteme1} by  \eqref{equ:fequaldiag}. Using the two first equations in item (1), the result becomes true on any ribbon graph.
 
\qed

\vspi 

It is natural to find the restricted polynomial $\cR'$ \eqref{rpr} over classes of HR-equivalent pinched ribbon graphs and 
to show for $\cR'$ corresponding recipe theorems (see Appendix \ref{app:recipinch}). 

\section{Topological transition polynomial}
\label{sect:trans}

We focus now on another interesting extension
of the so-called transition polynomial  for embedded
graphs \cite{joan,joan1} to HERGs. Note that, in this section, we will introduce
this transition polynomial at the level of HERGs and not at the level
of their HR-classes. Only at the end, the extension to 
classes will be quickly discussed.

\subsection{Vertex and HEGs states}
\label{subsect:vexstat}

\begin{figure}[htbp] 
\includegraphics[scale=0.8]{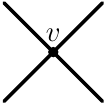}
\hspace{1cm}
\includegraphics[scale=0.8]{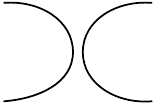}
\hspace{1cm}
   \includegraphics[scale=1.1]{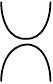}
\hspace{1cm}
   \includegraphics[scale=0.8]{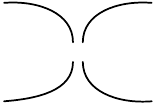}
  \caption{ Vertex states of a vertex $v$ (most left) in a HEG: 
(in that order) the black, white and cutting states of $v$.}
\label{figure: state}
\end{figure}

	Consider a $4$-regular graph with half-edges and $v$ one of its vertices. A vertex state at a vertex $v$ is a partition (with at most two elements in any subset of this partition) of the half-edges incident to $v$. Here, the vertex state is an extension of that given in \cite{joan} for graphs without half-edges since, in that case, the partition is strictly restricted into  pairs of half-edges incident to $v$. 
A subset in any such partition containing only one element will correspond to a half-edge which can be linked or not to
another half-edge in some other vertex state. 
Thus here, a vertex state here consists still in a set of disjoint curves.
For a given vertex in a  HEG, we obtain the usual three (white smoothing, black smoothing and crossing) states obtained when the graph is without half-edges \cite{joan} but also more states. These latter are obtained by cutting one or the two arcs given in the first three states considered. There is one particular state which will turn out to 
be important in the next analysis. We call it the
 cutting state and it is  obtained by cutting the two arcs 
of a black (or white) smoothing state (see Figure \ref{figure: state}).
The states that will retain our attention are the
black, white smoothing and cutting states.  
The rest of the states (even the crossing state) 
will not be used in the following. 
It might be  however interesting to  incorporate
those states in a more general study. 

For an abstract graph, the black and white smoothing 
(together with crossing) states cannot be distinguished. 
However, for an embedded graph,  choosing an appropriate
neighborhood of the half-edges in the partitioned vertex  
leads to different configurations.

Consider a  4-regular graph $F$ with half-edges,  
and a choose a particular
configuration of vertex states at each of its vertices. 
Then one calls this resulting graph configuration a graph state $S$
of $F$. 
We denote $c(S)$ the number of components of the 
graph state $S$ which are closed  and $c'(S)$ the number of those which are open. Importantly, in the following, the cutting state will replace 
the black smoothing state. As a mapping between these two
states,  we can  insert 
fictitious vertices of degree 2 between particular pairs of half-edges coming
from the cutting state, in such a way that these pairs mimic the pairs
of half-edges of the black smoothing. 
Now, $c''(S)$ counts the number of closed curves having at least one
fictitious vertex of degree 2. We will call these punctured curves. 
Note that, removing all the 2-valent vertices from punctured curves, 
one must end up with the open curves (in other words, open curves issued from the cutting state can be always glued together by 2-valent vertices and must form punctured curves).

	Let $F$ be a $4$-regular graph with half-edges.  A weight system $W(F)$ of $F$ is an assignment of weight valued in  a unitary ring $\mathfrak{R}$ to every vertex state of $F$. The state weight of a graph state $S$ of $F$ with weight system $W$ is $w(S)=\prod_vw(v,S)$
where $w(v,S)\in \mathfrak{R}$ is the vertex state weight of the vertex state at $v$ in the graph state $S$. 

\begin{definition}
Let $F$ be a $4$-regular graph with half-edges having weight system $W$ with values in the unitary ring $\mathfrak{R}$. Then the state model formulation of the generalized transition polynomial is given by
\beq \label{Def:transpoly}
q(F;W,a,b,d)=\sum_{S}w(S)a^{c(S)}b^{c'(S)}d^{c''(S)},
\eeq
where the sum is over all graph states $S$ of $F$.
\end{definition}

As an easily checked property, 
one proves that $q$ is multiplicative for
the disjoint union operations on 4-regular graphs. 

We now define the notion of medial graph with half-edges.
If $G$ is a HEG and it is cellularly embedded, its medial graph $G_{m}$ is constructed by 

- placing a vertex of degree 4 on each edge of $G$
and placing a ``fictitious'' vertex of degree 2 at the end point of each 
half-edge

- and drawing the edges  of the medial graph by following the face boundaries of $G$. 

- then remove all vertices of degree 2.

The insertion of vertices of degree 2 allows to keep track
of the open faces and new types of closed components
formed by these open faces. 
By convention, the medial graph of an isolated vertex is an isolated closed face. In the case of $n$ half-edges attached to this isolated vertex, the medial graph is a collection of open faces between with $n$ 2-valent vertices.  
	
\begin{figure}[htbp] 
\includegraphics[scale=2]{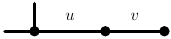}
\hspace{1cm}
   \includegraphics[scale=2]{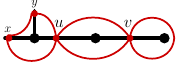}
\hspace{1cm}
   \includegraphics[scale=2]{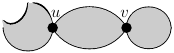}
  \caption{A HEG  (left), its medial graph construction in red (middle) 
with insertion of 2-valent vertices $x$ and $y$ which should
be removed at the end; its canonical checkerboard coloring (right).}
\label{figure: medial}
\end{figure}

A checkerboard coloring of a  cellularly embedded HEG is an assignment of the color black or white to each face  such that adjacent faces receive different colors.   We can color the faces of the medial graph with half-edges containing a vertex of the original graph $G$ black and the remaining white. We keep calling this, the canonical checkerboard coloring of $G_m$. In Figure \ref{figure: medial},  
we illustrate this construction.

 Let $G_m$ be a canonically checkerboard colored medial graph
and $v$ one of its vertices. 
$(G_m)_{wh(v)}$ is the  embedded graph with half-edges resulting from taking the white smoothing at the vertex $v$ of $G_m$, $(G_m)_{bl(v)}$  the embedded graph with half-edges obtained  from taking the black smoothing at the vertex $v$, and $(G_m)_{cut(v)}$ the embedded graphs with half-edges resulting from taking the cutting state at the vertex $v$.
 More embedded graphs from other states at $v$ could
be introduced but they will be not important for the
rest of the present work.  Precisely, there are parameters assigned to the other states but, in our next developments, 
these parameters are set at value 0.

The following proposition is immediate at this stage:
\begin{proposition} \label{medial:check}
Let $G$ be an embedded graph (with half-edges) with embedded, canonically checkerboard colored  medial graph $G_m$, and let $e$  be any edge of $G$, with $v_e$ the associated vertex in $G_m$. Then
\begin{enumerate}
\item $(G_m)_{bl(v_e)}=(G-e)_m$,
\item $(G_m)_{wh(v_e)}=(G/e)_m$,
\item $(G_m)_{cut(v_e)}=(G\vee e)_m$.
\end{enumerate}
\end{proposition}

  Consider an embedded graph $G$ with half-edges and its medial graph $G_m$. A vertex $v\in G_m$ has some state weight given by ordered elements $(\boldsymbol\alpha, \boldsymbol\eta)$, specifying the weight of the white smoothing state (referred as `uncut' vertex state),  
and cutting state, in that order.

\subsection{Topological transition polynomial for graphs with half-edges}

We have all  prerequisites to introduce a new transition polynomial.

\begin{definition} \label{def:trans}
 Let $G$ be an embedded graph with half-edges with canonically checkerboard colored embedded medial graph $G_m$, and let $W_m(G_m)=(\boldsymbol\alpha,\boldsymbol\eta)$ be a weight system associated with $G_m$. Then the \emph{topological transition polynomial} for  $G$ expresses as
\bea
Q(G, (\boldsymbol\alpha,\boldsymbol\eta), a, b,d) :=q(G_m; W_m,a,b,d).
\eea
\end{definition}

The following proposition is  straightforward. 

\begin{proposition}\label{p.recQ}
The topological transition polynomial may be computed by repeatedly applying  the following linear recursion relation at each vertex $v\in \cV(G_m)$, and, when there are no more vertices, evaluating each of the resulting closed curves 
to an independent variable $a$, the resulting open curves to an independent variable $b$, and the resulting punctured curves to 
an independent variable $d$: 
\beq\label{transip}
q(G_m, W_m,a,b,d)= \alpha_v q((G_m)_{wh(v)}, W_m,a,b,d)+ \eta_v q((G_m)_{cut(v)}, W_m,a,b,d).
\eeq
\end{proposition}

Note that, in practice during the calculation of this
recurrence rule and in order to be not confused
by open curves which might 
seem not be on the same connected component,
it is better to restore the fictitious 2-valent vertices.
Applying repeatedly \eqref{transip} when all the vertices of valence 4 have been decomposed, one removes then from the result all 2-valent vertices
and count the number of open and closed curves. 
We can provide an example as follows: Consider the
graph $\cG=$  
 $\includegraphics[angle=0, width=3.5cm, height=0.4cm]{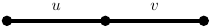}$, then $\cG_m=$ $\includegraphics[angle=0, width=3cm, height=1cm]{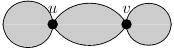}$. We apply \eqref{transip} successively as follows:  
\bea
Q(G, (\boldsymbol\alpha,\boldsymbol\eta), a, b.d)  &&= \alpha_{u} \includegraphics[angle=0, width=3cm, height=1cm]{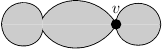} + 
\eta_{u} \includegraphics[angle=0, width=3cm, height=1cm]{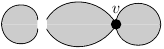}
\cr
&&= \alpha_{u} (\alpha_{v}\includegraphics[angle=0, width=3cm, height=1cm]{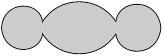} + \eta_{v} \includegraphics[angle=0, width=3cm, height=1cm]{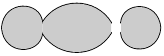}) +\cr 
&&\eta_{u} (\alpha_{v}\includegraphics[angle=0, width=3cm, height=1cm]{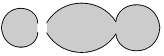} + \eta_{v} \includegraphics[angle=0, width=3cm, height=1cm]{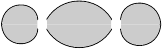}) \cr
&&= \alpha_{u}\alpha_{v} a + \alpha_{u}\eta_{v}b^{2}d^2 + \eta_{u}\alpha_{v}b^{2}d^2 + \eta_{u}\eta_{v}b^{4}d^3.
\eea

\vspi

\begin{theorem}
\label{theo:Qrec}
Let $G$ be an embedded graph with half-edges and $e\in \cE(G)$. Then 
\beq
Q(G;(\boldsymbol\alpha,\boldsymbol\eta),a,b,d)=\alpha_{e}Q(G/e;(\boldsymbol\alpha,\boldsymbol\eta),a,b,d) + \eta_{e}Q(G\vee e;(\boldsymbol\alpha,\boldsymbol\eta),a,b,d),\nonumber
\eeq
where, on the right-hand side, $(\boldsymbol\alpha,\boldsymbol\eta)$ denotes the weight system for $G$ restricted to $G/e$ or $G\vee e$ which is obtained by eliminating the weights for the vertex associated
with the edge $e$.
\end{theorem}
\proof
From Definition \ref{def:trans} and Propositions \ref{p.recQ} and \ref{medial:check}, the result follows.  
\qed
	
\begin{proposition}
In any ribbon graph with half-edges $G$ 
and $G_m$ its medial graph, we have:
\bea
F_{\ext}(G)=f(G) =c'(G_m), 
\qquad 
C_\partial(G) = c''(G_m), 
\eea
 where $F_{\ext}(G)$ is the number of external faces, $f(G)$ the number of half-edges and $C_\partial(G)$ the number of connected components of the boundary of $G$; 
 $c'(G_m)$ is the number of open
curves in $G_m$ and $c''(G_m)$ its number of punctured curves.
\end{proposition}
\proof
Notice that a connected component of the boundary is obtained by following a certain number of half-edges and external faces which alternate.
And to such a connected component we have a unique 
punctured curve. 

\qed

The following statement holds: 

\begin{theorem} \label{Theo:transpoly}
Let $\cG$ be a ribbon graph with half-edges with topological medial graph $\cG_m$. Then, fixing all $\alpha_v=\alpha$ and $\eta_v =\eta$, for all $v$,
we have
\bea
Q(\cG;(\boldsymbol\alpha,\boldsymbol\eta),a,b,d)=\alpha^{r(\cG)}\eta^{n(\cG)}a^{k(\cG)}\mathcal{R}(\cG;\frac{\eta a}{\alpha}+1,\frac{\alpha a}{\eta}+1,\frac{1}{a},d,b,1).
\eea
\end{theorem}
 \proof
Consider $\cG$ a ribbon graph possibly with half-edges and 
$\cG_m$ the  associated topological medial graph. Let $S$
be a state graph of $\cG_m$, 
$u(S)$ be its number of uncut vertex states (white-smoothing) and $q(S)$ its number of cut vertex states (cutting state). 
Introduce a weight system  $W_m=(\boldsymbol\alpha,\boldsymbol\eta)$ on $\cG_m$ such that $\alpha_v=\alpha$ and $\eta_v =\eta$, for all $v$, and define 
 $Q(\cG;(\boldsymbol\alpha,\boldsymbol\eta),a,b,d)= q(\cG_m;W_m,a,b,d)=\sum_{S\in \mathcal{S}}\alpha^{u(S)}\eta^{q(S)}a^{c(S)}b^{c'(S)}d^{c''(S)}$, with $\mathcal{S}$
the set of  graph states. Consider $A\sset \cG$, the graph state $S$ associated with $A$ is obtained by taking a cut vertex state for all vertices $v_e$ with $e\notin \cE(A)$ and uncut state for the remaining
ones, i.e. for $v_e$ such that $e \in \cE(A)$. For each uncut (resp. cut) vertex state,  we assign an uncut (resp. cut) edge in $\cG$. Hence $Q$ can be written as, given $\cE(\cG)=E$, $|\cE(A)|= |A|$, $|\cV(\cG)|=v(\cG)$, 
\bea
Q(\cG;(\boldsymbol\alpha,\boldsymbol\eta),a,b,d)&=&\sum_{A\sset \cG}\alpha^{|A|}\eta^{|E|-|A|}a^{F_{\inter}(A)}b^{F_{\ext}(A)} d^{C_\partial(A)}\crcr
&=&\eta^{|E|}\sum_{A\sset \cG}\Big(\frac{\alpha}{\eta}\Big)^{|A|}a^{F_{\inter}(A)} d^{C_{\partial}(A)} b^{F_{\ext}(A)}.
\eea 
On the other hand, 
 \bea
&&
\mathcal{R}(\cG;x,y,z,s,t,1)
=(x-1)^{-k(\cG)}(y-1)^{-v(\cG)}z^{-v(G)}\times 
\cr\cr
&&\sum_{A\sset \cG}((x-1)(y-1)z^2)^{k(A)} ((y-1)z)^{E(A)}z^{-F_{\inter}(A)}s^{C_\partial(A)}t^{f(A)}.
\eea
 Setting $z=\frac{1}{a}$, $y=\frac{\alpha a}{\eta}+1$, $x=\frac{\eta a}{\alpha}+1$, $s=d$ and $t=b$, one has: 
\bea
\mathcal{R}(\cG;x,y,z,s,t,1)&=&\Big(\frac{\eta}{\alpha}\Big)^{v(\cG)-k(\cG)}a^{-k(\cG)}\sum_{A\sset\cG}\Big(\frac{\alpha}{\eta}\Big)^{E(A)}a^{F_{\inter}(A)} d^{C_\partial(A)} b^{f(A)}\cr
&=&(\eta)^{-n(\cG)}(\alpha)^{-\rk(\cG)}a^{-k(\cG)}Q(\cG;(\boldsymbol\alpha,\boldsymbol\eta),a,b,d).
\eea

\qed

Let us now discuss how the above properties extend to HR-classes. 
This can be understood in a direct way: 
we can set a definition $Q([\cG]; -):=Q(\cG;-) = \cR(\cG) = \cR([\cG])$,
and $Q([\cG;-])= q([\cG_m], -)$  where we define the class $[\cG_m]$ 
of medial graphs obtained by choosing $\cG_m$ the medial graph  
of $\cG$ and all medial graphs obtained from $\cG_m$ by displacing the fictitious 2-valent vertices within the same punctured curve
or from one punctured curve to the other. 
We must keep however at least one fictitious 2-valent
vertex per punctured curve. A moment of thought 
leads to the correspondence between these latter moves 
with the HR-moves. This can be
used to achieve the claim $\forall g \in [\cG_m]$, 
$q(g;-) = q(\cG_m;-)$. Thus we can define $q([\cG_m], -)=q(\cG_m;-)
= Q(\cG;-)= Q([\cG];-)$. 

\vspi

Several other interesting developments can be now undertaken from the polynomial invariants treated in this paper. For instance, 
we might ask if one can make sense of a duality relation for the generalized transition polynomial defined on HERGs. Finding a dual for a HERG becomes however a nontrivial task. Indeed, two HERGs belonging
to the same HR-equivalence class do not always have equivalent dual. 
These certainly deserves to be elucidated. On another connected
domain, significant progresses around matroids  \cite{Duchamp:2013joa} and  Hopf algebra techniques \cite{Duchamp:2013pha} applied to the Tutte polynomial have been recently highlighted. These studies should find as well an extension for the types of invariants worked out in
this paper. Finally, combining some ideas of the present work and Hopf algebra calculations as found in \cite{Raasakka:2013kaa}, 
one might be able to prove a universality theorem for polynomial invariants over stranded graphs \cite{rca} generalizing ribbon graphs with 
half-edges.

\vspi 

\begin{center}
{\bf Acknowledgements}
\end{center}

 This work is partially supported by the Abdus Salam International Centre for Theoretical Physics (ICTP, Trieste, Italy) through the Office of External Activities (OEA)-Prj-15. The ICMPA is also in partnership with the Daniel Iagolnitzer Foundation (DIF), France.
JBG acknowledges the support of the Alexander von Humboldt Foundation. 

\section*{ Appendix}
\label{app}

\appendix

\renewcommand{\theequation}{\Alph{section}.\arabic{equation}}
\setcounter{equation}{0}

\section{Recipe theorems for the polynomial on pinched ribbon graphs}
\label{app:recipinch}

  We provide here two recipe theorems for the
 polynomial $\cR'$ defined HR-classes of pinched ribbon graphs. The 
proof of the first recipe theorem is very similar to the proof as found in \cite{joan}. This is mainly due to the factorization property of $\cR'$ with respect to the one-point-joint operation.
We also use particular diagrams called canonical chord-diagrams with half-chords associated with  the pinched one-vertex ribbon graphs with half-edges. Notice that half-chords here can be disposed in arbitrary way on the circle. The  variable $t$ appears as a  deformation parameter
during this proof. 
The second recipe holds for a more general 
contraction/deletion rule weighted by parameters. It can be considered 
as the most  complete recipe theorem in that matter for $\cR'$. 
For completeness reasons, we report 
these theorems for the interested reader. 

The change of variable $S\to Z^{-1}$ in $\cR$
leads to the polynomial $\cR'$ defined by:
\beq
\cR_{\cG}(X,Y,Z,Z^{-1},S,W,T) = \cR'_{\cG}(X,Y,Z,W,T) \,.
\eeq
A quick look of the expression of $\cR'$, one infers the relation
\beq
 \cR'_{\cG}(X,Y,Z,W,T)=T^{|\mf^0|}T^{2n(\cG)} \,R_{\cG}(\tilde X, \frac{\tilde Y}{T^2},Z,W)\,,
\eeq 
where
\begin{eqnarray}
\left\{\begin{array}{ll}\tilde X=(X-1)T^2+1,\\
 \tilde Y=Y-1+T^2. \end{array}
\right.
\end{eqnarray}

As a polynomial on a different category of 
graphs, a recipe theorem for $\cR'$ on
pinched ribbon graphs can be investigated on its own right. 

To proceed with, let us define now some canonical chord diagrams with half-chords associated  to the HR-classes of pinched one-vertex ribbon graphs. Let  $D_{ijkl}$ be the following chord diagram called canonical: 
it consists of $i-2j-k$ positive disjoint chords,
$j$ pairs of intersecting positive chords,  $k$  number of   negative disjoint chords and  $l$  number of half-chords. Then, we have: 
\bea
\mathcal{R}'(D_{1000};x,y,z,w,t) = (y-1) + t^2,\\
\mathcal{R}'(D_{10010};x,y,z,w,t) = (y-1)zw + t^2,\\
\mathcal{R}'(D_{2100};x,y,z,w,t) = (y-1)^2z^2 + 2(y-1)t^2 + t^4,\\
\mathcal{R}'(D_{0001};x,y,z,w,t) = t.
\eea

Using the fact that the polynomial $\mathcal{R}'$ is multiplicative 
 for the one-point-joint of pinched ribbon graphs, we can write
\bea
\mathcal{R}'(D_{ijkl})=[\mathcal{R}'(D_{1000})]^{i-2j-k}[\mathcal{R}'(D_{2100})]^{j}[\mathcal{R}'(D_{1010})]^{k}[\mathcal{R}'(D_{0001})]^l.
\eea
\begin{theorem}[Recipe theorem for $\mathcal{R}'$]
\label{theo:Recipe1}
Let $\mathcal{M}$ be a minor closed subset of pinched ribbon graphs containing all pinched ribbon graphs on two vertices. Let $F$ be a map from $\mathcal{M}$ to a commutative ring $\mathcal{R}$ with unity.  Let $s=F(D_{2100})$, $q=F(D_{1000})$, $r=F(D_{1010})$ and  $s_1=F(D_{0001})$  and suppose there exist elements $\alpha$, $x$, $u$, $v$ $\in \mathfrak{R}$ with $\alpha$ a unit such that:
\begin{enumerate}
\item \bea \label{systempinch}
F(G) = \left\{\begin{array}{ll} 
F(G \vee e) + F(G/e) & {\text{if e is  regular}},\\\\
((x-1)(\alpha^{-1}s_1)^2 +1)  \, F(G/e) & {\text{if e is a bridge}}.
\end{array} \right. 
\eea
\item $F(G\sqcup H)=F(G)F(H)$ and $\alpha F(G\cdot H)=F(G)F(H)$ where 
$G$ and $H$ are pinched embedded bouquets.
\item $F(E)=\alpha^nt^m$ if $E$ is a pinched edgeless graph with $m$ 
half-edges and $n$ vertices;
\item $(q- s_1^2)^2u^2=\alpha[s-2\alpha^{-2}(q- s_1^2)s_1^2-  s_1^4]$, $(q-s_1^2)uv=r- s_1^2$, and $v=v^2$.
Then
\bea
F(G)=\alpha^{k(G)}\mathcal{R}'(G;x,\alpha^{-1}
(q-s_1^2+\alpha ),u,v,\alpha^{-1}s_1),\label{Recipe1}
\eea
where $k(G)$ is the number of components of $G$.
\end{enumerate}
\end{theorem}
\proof
The proof  proceeds in two steps: (i) we focus on pinched one-vertex ribbon graphs  and then (ii) on general pinched ribbon graphs.
\begin{itemize}
\item[(i)]
	Using items (2) and (4) of Theorem \ref{theo:Recipe1}, the relation \eqref{Recipe1} holds in the canonical 
diagrams with half-chords. Indeed, considering a canonical diagram $D_{ijkl}$ and using the relations in item (4), Theorem \ref{theo:Recipe1} is immediately verified for $D_{1000}$, $D_{2100}$, $D_{1010}$ and $D_{0001}$. From item (2), we have: 
\bea
\alpha^{i-j+l-1}F(D_{ijkl})=[F(D_{1000})]^{i-2j-k}[F(D_{2100})]^{j}[F(D_{1010})]^{k}[F(D_{0001})]^l.\label{equation:onepointjoint}
\eea
Using the relation \eqref{equation:onepointjoint} and the fact that the theorem holds for the 
diagrams $D_{1000}$, $D_{2100}$, $D_{1010}$ and $D_{0001}$, then
 it is  also true on any canonical diagram $D_{ijkl}$. From item (1), we obtain the relations:
\bea
F(D_1) - \mu F(D'_1) = F(D_2) - \mu F(D'_2)
 \label{Related1} \mbox{ and },\\
F(D_3) - \mu F(D'_3) = F(D_4) - \mu F(D'_4),
 \label{Related2}
\eea
where  $D_1$, $D_2$, $D_3$ and $D_4$ are related as shown in Figures \ref{oprotate} and \ref{optwiste} \footnote{Notice here that, the position of the half-chords does not matter. We can also use diagrams without half-chords for illustration.}
with $D'_i=D_i\vee e$, and $\mu=1$ if there is a chord from $a\cup b$ to 
$c\cup d$,  otherwise $\mu=[(x-1)(\alpha^{-1}s_1)^2 +1]$.
 The relation \eqref{Related1} holds also for $\mathcal{R}'$ and hence for
 $F'=F - \alpha \mathcal{R}'$.

	By induction, assume that $F(D) = \alpha \mathcal{R}'(D)$ for any signed chord diagram with fewer than $n$ chords. Then $F'$ vanishes on any signed chord diagram with fewer than $n$ chords.
	From  \eqref{Related1}, we have $F'(D)=F'(D_{njkl})$ where $D$ is related to a canonical diagram $D_{njkl}$. From the following equations:
\bea
F(D_{njkl})=\alpha\mathcal{R}'(D_{njkl};x,\alpha^{-1}
(q-s_1^2+\alpha ),u,v,\alpha^{-1}s_1)
\eea
obtained using \eqref{equation:onepointjoint}, and
\bea
F'(D)=F'(D_{njkl})=0,
\eea
	we get, by induction, $F'(D)=0$ on any signed chord diagram $D$.
 Finally, the result holds on any rosette of pinched ribbon graph. 
\item[(ii)]Using item (1), the result becomes obvious on any pinched ribbon graph.

\qed

\end{itemize}

In this proof we remark that, we do not need $\mathcal{M}$ to contain all ribbon graphs on two vertices. In fact, it is sufficient to consider that  $F$  satisfies  \eqref{Related1} and \eqref{Related2} for a chain of rosettes terminating in a canonical chord diagram.
The next definition will be needed. 

We say that a minor closed set $\mathcal{M}$ of pinched ribbon graphs is closed under chord operations whenever $D\in\mathcal{M}$ and $D\sim D_{ijkl},$ then there is a finite sequence $D=D_1...D_n=D_{ijkl}$ with $D_i\in\mathcal{M}$ and $D\sim D_{i+1}$ for all $i$.

\begin{theorem}[The ``low fat'' recipe theorem for $\cR'$]
\label{theo:lowrecipeprime} 
Theorem \ref{theo:Recipe1} holds with ``Let $\mathcal{M}$ be a minor closed subset of pinched ribbon graphs containing all ribbon graphs on two vertices and let $F$ map $\mathcal{M}$ to a commutative ring $\mathfrak{R}$ with unity'', replaced by ``Let $\mathcal{M}$ be a minor closed subset of pinched ribbon graphs closed under chord operations that contains $D_{1000}$, $D_{1010}$, $D_{2100}$, $D_{0001}$ and let $F$ map $\mathcal{M}$ to a commutative ring $\mathfrak{R}$ with unity, such that $F$ satisfies  \eqref{Related1} and \eqref{Related2} whenever the $D_i$'s are related as in Figures \ref{oprotate} 
and \ref{optwiste}.''
\end{theorem}

From this point, we also have the following statement: 
\begin{corollary}
If $F$, $\mathcal{M}$, $\mathfrak{R}$  satisfy the conditions of Theorem \ref{theo:lowrecipeprime}, with both $
q-s_1^2$ and $r-s_1^2$ being unit of $\mathfrak{R}$, then $v=1$, and thus $F$ does not discern orientation by the presence or absence of a single idempotent element.
\end{corollary}
\proof From the second equation in the above item (4), we infer
$\frac{r- s_1^2}{q- s_1^2}=uv=
uv^2=v\frac{r-s_1^2}{q-s_1^2}.$ 
Furthermore, $v=1$ since $q-s_1^2$ and 
$r-s_1^2$ are units. Therefore,
$F(G)=\alpha^{k(G)}\mathcal{R}'(G;x,\alpha^{-1}
(q-s_1^2+\alpha ),u,1,\alpha^{-1}s_1)$.

\qed

 There exists another statement about the recipe theorem
when the contraction/deletion rule is weighted by parameters. 
We now provide, in that more general case, a second more generic
recipe to construct a function satisfying the weighted recurrence
relation. 

\begin{theorem}[Second recipe theorem for $\mathcal{R}'$] 
\label{theo:recipe1}
Let $\mathcal{M}$ be a minor closed subset of pinched ribbon graphs containing all pinched ribbon graphs on two vertices. Let $F$ be a map from $\mathcal{M}$ to a commutative ring $\mathcal{R}$ with unity.  Let $s=F(D_{2100})$, $q=F(D_{1000})$ and $r=F(D_{1010})$; and suppose there exist elements $\alpha$, $x$, $u$, $v$ $\in \mathfrak{R}$ with $\alpha$ a unit such that:
\begin{enumerate}
\item \bea \label{systempinch2}
F(G) = \left\{\begin{array}{ll} 
\sigma F(G-e) + \tau F(G/e) & {\text{if e is  regular}},\\\\
x \, F(G/e) & {\text{if e is a bridge}}.
\end{array} \right. 
\eea
\item $F(G\sqcup H)=F(G)F(H)$ and $\alpha F(G \cdot H)=F(G)F(H)$ where $G$ and 
$H$ are pinched embedded bouquets.
\item $F(E)=\alpha^nT^m$ if $E$ is an edgeless graph with $m$ half-edges and $n$ vertices.
\item $(q-\alpha\tau)^2u^2=\alpha(s-2\tau q+\alpha\tau^{2})$, $(q-\alpha\tau)uv=r-\alpha\tau$, and $v=v^2$.
Then
\bea
&& 
\sigma^{-r(G)}\tau^{-n(G)}F(G)=\alpha^{k(G)}R(G;\sigma^{-1}x,\alpha^{-1}\tau^{-1}q,u,v)\cr\cr
&&= T^{-(|\mf^0|+2n(\cG))}\alpha^{k(G)}\mathcal{R}'(\sigma^{-1}[(x-\sigma)T^{-2}+ \sigma],T^{2}(\alpha^{-1}\tau^{-1}q)-T^{2}+1,u,v,T),\label{recipe1}
\eea
where $k(G)$, $r(G)$ and $n(G)$ are respectively the number of connected components, rank and nullity of $G$.
\end{enumerate}
\end{theorem}
\proof This theorem holds on any canonical diagram $D_{ijkl}$ or $D_{ijk}$ (BR diagrams).
 Using item (1), we have the relation
\bea
\tau F(D_1) - \mu\sigma F(D'_1) = \tau F(D_2) - \mu\sigma F(D'_2)
 \label{rlated1} \mbox{ and }\\
\tau F(D_3) - \mu\sigma F(D'_3) = \tau F(D_4) - \mu\sigma F(D'_4)
 \label{rlated2},
\eea
where  $D_1$, $D_2$, $D_3$ and $D_4$ are related as shown in Figures \ref{oprotate} and \ref{optwiste} with $D'_i=D_i \vee e$, and $\mu=1$ if there is a chord from $a\cup b$ to $c\cup d$, and otherwise $\mu=x$. The polynomials  $\mathcal{R}'$ (respectively $R$), satisfies  \eqref{rlated1} and \eqref{rlated2} and hence 
 $F'(G)=\sigma^{-r(G)}\tau^{-n(G)}F(G) - T^{-(|\mf^0|+2n(\cG))}\alpha^{k(G)} \mathcal{R}'(G)$ (respectively $F'(G)=\sigma^{-r(G)}\tau^{-n(G)}F(G) - \alpha R(G)$) satisfies the same relation.
Assume that $\sigma^{-r(D)}\tau^{-n(D)}F(D) = T^{-(|\mf^0|+2n(D))}\alpha^{k(D)} \mathcal{R}'(D)$ (respectively $\sigma^{-r(D)}\tau^{-n(D)}F(D) = \alpha R(D)$)  for any signed chord diagram with fewer than $n$ chords. Therefore $F'$ vanishes on any signed chord diagram with fewer than $n$ chords.
	Using  \eqref{rlated1} and \eqref{rlated2} we have $\tau F'(D)=\tau F'(D_{njkl})$ (respectively $\tau F'(D)=\tau F'(D_{njk})$) where $D$ is related to a canonical diagram $D_{njkl}$ (respectively $D_{njk}$). From the following equations:
\bea
\tau^{-n}F(D_{njkl})=T^{-l-2n}\alpha\mathcal{R}'(D_{njkl})=\alpha R(D_{njk}),
\eea
and
\bea
F'(D)=F'(D_{njkl})=0,
\eea
we have by induction, $F'(D)=0$ on any signed chord diagram $D$.
 Finally the result holds on any rosette of pinched ribbon graph. 
Using item (1), the result becomes true on any pinched ribbon graph.

\qed

We can introduce the weaker form of the recipe theorem.

 \begin{theorem}[The ``low fat'' recipe theorem for $\cR'$]
\label{theo:lowrecipeprime1} 
Theorem \ref{theo:recipe1} holds with ``Let $\mathcal{M}$ be a minor closed subset of pinched ribbon graphs containing all ribbon graphs on two vertices and let $F$ map $\mathcal{M}$ to a commutative ring $\mathfrak{R}$ with unity'', replaced by ``Let $\mathcal{M}$ be a minor closed subset of pinched ribbon graphs closed under chord operations that contains $D_{1000}$, $D_{1010}$, $D_{2100}$ and let $F$ map $\mathcal{M}$ to a commutative ring $\mathfrak{R}$ with unity, such that $F$ satisfies  \eqref{rlated1} and \eqref{rlated2} whenever the $D_i's$ are related as in Figures \ref{oprotate} and \ref{optwiste}.''
\end{theorem} 

\begin{corollary}
If $F$, $\mathcal{M}$, $\mathfrak{R}$ satisfy the conditions of Theorem \ref{theo:lowrecipeprime1}, with both $q-\alpha \tau$ and $r-\alpha\tau$ being unit of $\mathfrak{R}$, then $v=1$, and thus $F$ does not discern orientation by the presence or absence of a single idempotent element.
\end{corollary}

\proof
From the second equation in the above item (4), we infer
$\frac{r-\alpha\tau}{q-\alpha\tau}=uv=
uv^2=v\frac{r-\alpha\tau}{q-\alpha\tau}.$ 
Furthermore, $v=1$ since $q-\alpha\tau$ and 
$r-\alpha\tau$ are units. Therefore,
$\sigma^{-r(G)}\tau^{-n(G)}F(G)=\alpha^{k(G)}\mathcal{R}'(G;\sigma^{-1}x,\alpha^{-1}\tau^{-1}q,u,1)$.

\qed

\section{One-point-joint and semi-factorization property of $\cR$}
\label{app:semifacto}

In this appendix, we investigate the precise breaking of the factorization 
property of $\cR$ for the one-point-joint operation
of two HERGs or of two HR-classes. First, a precision must be given on the notation 
$\cG_1 \cdot_{v_1,v_2} \cG_2$ and we will adopt instead
the following notation:
\bea
\cG_1 \cdot_{(v_i,c_i)} \cG_2\,, 
\eea
where the pairs $(v_i,c_i)$, $i=1,2$, incorporate the vertex $v_i$ of $\cG_i$ where the joining operation is performed and an arc
 $c_i$ on $v_i$ (which without edges)  where 
$v_1$ and $v_2$ merge.

Let us denote $\cP_{\cG}$ the set of cutting spanning subgraphs of $\cG$. 
We denote $\mh_{A_i}(c_i)$ the face of $A_i \sset \cG_i$ which 
contains the arc $c_i$ in $v_i$ seen as a vertex of $A_i$. 
And let us call $\cP_{\cG_i}(v_i,c_i,\eps)$, $\eps\in \{ 0,1\}$, 
the set of cutting spanning subgraphs of $\cG_i$ 
such that 
\bea
&&
A_i \in \cP_{\cG_i}(v_i,c_i,0)\subset \cP_{\cG_i} \qquad \Leftrightarrow
\qquad \mh_{A_i}(c_i) \in \cF_{\inter}(A_i) \,,\crcr
&&
A_i \in \cP_{\cG_i}(v_i,c_i,1)\subset \cP_{\cG_i} \qquad \Leftrightarrow
\qquad \mh_{A_i}(c_i) \in \cF_{\ext}(A_i) \,.
\eea

Consider a cutting spanning subgraph $A_i$ of $\cG_i$, 
when we write $A_1 \cdot_{(v_i,c_i)} A_2$, we implicitly refer
to the fact that $c_i$ must be seen as an arc of $v_i$ in $A_i$
(and not anymore as an arc of $v_i$ in $\cG_i$, except when,
of course, $A_i = \cG_i$).

One must see that  
the set of cutting spanning subgraphs of $\cG_1 \cdot_{(v_i,c_i)} \cG_2$
is one-to-one with the set of cutting spanning subgraphs of
$\cG_1 \sqcup \cG_2$, i.e. to each
$ A_1 \cdot_{(v_i,c_i)} A_2 $ corresponds a unique $A_1 \sqcup A_2$,
for $A_i \sset \cG_i$. It follows the  decomposition lemma:

 \begin{lemma}[Cutting spanning subgraph decomposition]
\label{cutspa}
Let $\cG_i$ two disjoint graphs, then the set $\cP_{\cG_1 \cdot_{(v_i,c_i)} \cG_2}$
is the union of two disjoint sets  $\cP_{\cG_1 \cdot_{(v_i,c_i)} \cG_2}^{0}$
and $\cP_{\cG_1 \cdot_{(v_i,c_i)} \cG_2}^{1}$ such that
\bea
\cP_{\cG_1 \cdot_{(v_i,c_i)} \cG_2}^{0} &\equiv& 
(\cP_{\cG_1}(v_1,c_1,0)\cup \cP_{\cG_2}(v_2,c_2,0))
\cup 
(\cP_{\cG_1}(v_1,c_1,0)\cup \cP_{\cG_2}(v_2,c_2,1))\cr\cr
&& 
\cup 
(\cP_{\cG_1}(v_1,c_1,1)\cup \cP_{\cG_2}(v_2,c_2,0)) \,,\cr\cr
\cP_{\cG_1 \cdot_{(v_i,c_i)} \cG_2}^{1} &\equiv& 
\cP_{\cG_1}(v_1,c_1,1)\cup \cP_{\cG_2}(v_2,c_2,1)\,.
\eea
where $\equiv$ means ``one-to-one with''. 
\end{lemma}
\begin{lemma}\label{ranktoface}
Let $A_i \sset \cG_i$, then 
\bea
&&
r(A_1 \cdot_{(v_i,c_i)} A_2) = r(A_1) + r(A_2) \,,\crcr
&&
n(A_1 \cdot_{(v_i,c_i)} A_2) = n(A_1) + n(A_2)\,, \crcr
&&
t(A_1 \cdot_{(v_i,c_i)} A_2) = t(A_1) + t(A_2) \,,\crcr
&&
k(A_1 \cdot_{(v_i,c_i)} A_2) = k(A_1) + k(A_2) -1\,,
\label{rntk}\\
&&
\hspace{-0.2cm} F_{\inter}(A_1 \cdot_{(v_i,c_i)} A_2) =
\left\{\begin{array}{l}
F_{\inter}(A_1) + F_{\inter}(A_2)-1 \,, \quad \text{if}\;\; A_1\cdot_{(v_i,c_i)} A_2  \in \cP_{\cG_1 \cdot_{(v_i,c_i)} \cG_2}^{0} \\
F_{\inter}(A_1) + F_{\inter}(A_2) \,, \quad \text{if}\;\;  A_1\cdot_{(v_i,c_i)} A_2  \in \cP_{\cG_1 \cdot_{(v_i,c_i)} \cG_2}^{1} 
\end{array} \right. \cr\cr
&& 
\hspace{-0.2cm}
C_{\partial}(A_1 \cdot_{(v_i,c_i)} A_2) =
\left\{\begin{array}{l}
C_{\partial}(A_1) + C_{\partial}(A_2) \,, \quad \text{if}\;\; A_1\cdot_{(v_i,c_i)} A_2  \in \cP_{\cG_1 \cdot_{(v_i,c_i)} \cG_2}^{0} \\
C_{\partial}(A_1) +C_{\partial}(A_2) -1 \,, \quad \text{if}\;\;  A_1\cdot_{(v_i,c_i)} A_2  \in \cP_{\cG_1 \cdot_{(v_i,c_i)} \cG_2}^{1} 
\end{array} \right. 
\eea
\end{lemma}
\proof The relations involving the rank, the nullity, the number of half-edges 
and number of connected components in \eqref{rntk} are directly obtained
from the usual additivity of these quantities during the one-point-joint. 

The relation involving the number of internal lines 
and connected components of the boundary can be  
inferred as follows: in the sector
$\cP_{\cG_1 \cdot_{(v_i,c_i)} \cG_2}^{1}$, after the one-point-joint,
one always looses a connected
component of the boundary. The contact
between $A_1$ and $A_2$ happens on two external faces
belonging necessarily to two initial and different connected components
 of the boundary of $A_1$ and $A_2$. Note that 
closed faces are all preserved. For the sector $\cP_{\cG_1 \cdot_{(v_i,c_i)} \cG_2}^{0}$, the reasoning is similar: one looses
an internal face this time for
each subsector  $\cP_{\cG_1}(v_1,c_1,0)\cup \cP_{\cG_2}(v_2,c_2,0)$, $\cP_{\cG_1}(v_1,c_1,0)\cup \cP_{\cG_2}(v_2,c_2,1)$ and $\cP_{\cG_1}(v_1,c_1,1)\cup \cP_{\cG_2}(v_2,c_2,0)$.
In the end, there is always an internal face which is lost
but the number of connected components of the boundary
is always preserved for each sector. 

\qed  

We are in position to understand how the BR polynomial for HERGs behaves under one-point-joint. Let us adopt a slightly different notation for the same object $\cR$ making explicit its dependence on the set of its cutting subgraphs
\bea
\cR(\cG;\cP_{\cG} ;x,y,z,s,w,t) = \sum_{A\in \cP_{\cG}} M_{\cR}(A; x,y,z,s,w,t)
\eea
where $M_{\cR}$ is the monomial ordinarily associated
with each subgraph and defined by $\cR$. 
Now, let us compute using Lemmas \ref{cutspa} and \ref{ranktoface}
\bea
&&
\cR(\cG_1 \cdot_{(v_i,c_i)} \cG_2;\cP_{\cG_1 \cdot_{(v_i,c_i)} \cG_2} )
 = 
\sum_{A_1 \cdot_{(v_i,c_i)}  A_2 \in \cP_{\cG_1 \cdot_{(v_i,c_i)} \cG_2}^{0} } 
M_{\cR}(A_1 \cdot_{(v_i,c_i)}  A_2; x,y,z,s,w,t) \crcr
&& 
 + 
\sum_{A_1 \cdot_{(v_i,c_i)}  A_2 \in \cP_{\cG_1 \cdot_{(v_i,c_i)} \cG_2}^{1} } 
M_{\cR}(A_1 \cdot_{(v_i,c_i)}  A_2; x,y,z,s,w,t) \cr\cr\cr\cr
&&
 = 
\sum_{A_1 \cdot_{(v_i,c_i)}  A_2 \in \cP_{\cG_1 \cdot_{(v_i,c_i)} \cG_2}^{0} } 
M_{\cR}(A_1; x,y,z,s,w,t)
M_{\cR}( A_2; x,y,z,s,w,t) \crcr
&& 
 + (zs)^{-1}
\sum_{A_1 \cdot_{(v_i,c_i)}  A_2 \in \cP_{\cG_1 \cdot_{(v_i,c_i)} \cG_2}^{1} } 
M_{\cR}(A_1; x,y,z,s,w,t)
 M_{\cR}(A_2; x,y,z,s,w,t) \crcr
&&
 = \Big\{ 
\sum_{A_1 \in \cP_{\cG_1}(v_1,c_1;0) } 
M_{\cR}(A_1; x,y,z,s,w,t) \Big\} 
\Big\{ \sum_{A_2 \in \cP_{\cG_2}(v_2,c_2;0) }  M_{\cR}( A_2; x,y,z,s,w,t) \Big\} \crcr
&& 
+ 
\Big\{ 
\sum_{A_1 \in \cP_{\cG_1}(v_1,c_1;0) } 
M_{\cR}(A_1; x,y,z,s,w,t) \Big\} 
\Big\{ \sum_{A_2 \in \cP_{\cG_2}(v_2,c_2;1) }  M_{\cR}( A_2; x,y,z,s,w,t) \Big\} \crcr
&& 
+ 
\Big\{ 
\sum_{A_1 \in \cP_{\cG_1}(v_1,c_1;1) } 
M_{\cR}(A_1; x,y,z,s,w,t) \Big\} 
\Big\{ \sum_{A_2 \in \cP_{\cG_2}(v_2,c_2;0) }  M_{\cR}( A_2; x,y,z,s,w,t) \Big\} \crcr
&&
 + (zs)^{-1}
\Big\{ 
\sum_{A_1 \in \cP_{\cG_1}(v_1,c_1;1) } 
M_{\cR}(A_1; x,y,z,s,w,t) \Big\} 
\Big\{ \sum_{A_2 \in \cP_{\cG_2}(v_2,c_2;1) }  M_{\cR}( A_2; x,y,z,s,w,t) \Big\}.
\cr
&&
\eea
This computes to give, in suggestive notation where $\cR(\cG, \cP)$
means that we compute the sum of monomials only for subgraphs in $\cP$, 
\bea 
&&
\cR(\cG_1 \cdot_{(v_i,c_i)} \cG_2;\cP_{\cG_1 \cdot_{(v_i,c_i)} \cG_2} )
 = \cr\cr
&& \cR(\cG_1;\cP_{\cG_1}(v_1,c_1;0)) \cR(\cG_2;\cP_{\cG_2}(v_2,c_2;0)) 
 +  \cR(\cG_1;\cP_{\cG_1}(v_1,c_1;1)) \cR(\cG_2;\cP_{\cG_2}(v_2,c_2;0))
\cr\cr
&&
+
 \cR(\cG_1;\cP_{\cG_1}(v_1,c_1;0)) \cR(\cG_2;\cP_{\cG_2}(v_2,c_2;1))
+
(zs)^{-1} \cR(\cG_1;\cP_{\cG_1}(v_1,c_1;1)) \cR(\cG_2;\cP_{\cG_2}(v_2,c_2;1)) 
\,.\crcr
&&
\eea
Thus, we obtain a broken factorization property of the
polynomial $\cR$ of the form:
\bea
&&
\cR(\cG_1 \cdot_{(v_i,c_i)} \cG_2)
 =  \cR(\cG_1) \, \cR(\cG_2;\cP_{\cG_2}(v_2,c_2;0)) 
\cr\cr
&&
+
 \Big(\cR(\cG_1;\cP_{\cG_1}(v_1,c_1;0)) +
(zs)^{-1} \cR(\cG_1;\cP_{\cG_1}(v_1,c_1;1)) \Big) \cR(\cG_2;\cP_{\cG_2}(v_2,c_2;1)) \cr\cr
&& 
 = \cR(\cG_1) \, \cR(\cG_2;\cP_{\cG_2}(v_2,c_2;0))
 + \widetilde\cR_{sz}(\cG_1;v_1,c_1) \cR(\cG_2;\cP_{\cG_2}(v_2,c_2;1))\,, 
\eea
where $\widetilde\cR_{sz}(\cG;v,c):= \cR(\cG;\cP_{\cG}(v,c;0)) +
(zs)^{-1} \cR(\cG;\cP_{\cG}(v,c;1)) $ is a deformed version
 of the $\cR(\cG)$ polynomial depending on a given vertex $v$ and
an arc $c$ on it. 
There exists of course a symmetric relation
where the role of $\cG_1$ is played by $\cG_2$. 
Finally, we do not see how the polynomial $\cR_{\cG_1\cdot \cG_2}$
 factorizes unless  $s=  z^{-1}$ which maps $\cR$
to $\cR'$. Finally, all the above properties extends to 
HR-classes without ambiguity.

\vspace{0.3cm}

\begin{center}
\rule{3cm}{0.01cm}
\end{center}

\end{document}